\documentclass[11pt]{amsart}

\usepackage{amsmath,amssymb,hyperref,color,enumerate,chngpage}
\IfFileExists{mathabx.sty}{\usepackage{mathabx}}{}
\usepackage{graphicx}

\usepackage{mathrsfs}
\usepackage{amsthm}

\newcommand{\beq} {\begin{equation}}
\newcommand{\eeq} {\end{equation}}
\newcommand{\bdm} {\begin{displaymath}}
\newcommand{\edm} {\end{displaymath}}

\newcommand{\bit}{\begin{itemize}}
\newcommand{\eit}{\end{itemize}}
\newcommand{\bde}{\begin{description}}
\newcommand{\ede}{\end{description}}

\newcommand{\ben}{\begin{enumerate}}
\newcommand{\een}{\end{enumerate}}

\newcommand{\barr}{\begin{array}}
\newcommand{\earr}{\end{array}}

\newcommand{\LRp}[1]{\left( #1 \right)}

\newcommand{\dive}{\operatorname{div}}

\newtheorem{theorem}{Theorem}[section]
\newtheorem{proposition}[theorem]{Proposition}
\newtheorem{lemma}[theorem]{Lemma}
\newtheorem{cor}[theorem]{Corollary}
\newtheorem{definition}[theorem]{Definition}
\newtheorem{rmk}[theorem]{Remark}

\usepackage{doi}

\title[Locking-free methods for poroelasticity optimal control problems]{Locking-free numerical methods for poroelasticity optimal control problems with minimal regularities}
\author{Jeonghun J. Lee}

\begin{document}

\begin{abstract}
We study distributed optimal control of a linear poroelastic system in a
total pressure formulation, with displacement and fluid-pressure tracking
and box constraints on the source control. A mechanical Schur complement
gives well-posedness at zero storage for each fixed finite Lam\'e parameter,
without estimates that grow like the inverse storage coefficient.
We discretize the state and adjoint with conforming finite elements and
dG(0) time stepping, and use piecewise constant controls.
Under the stated regularity assumptions and a positive lower bound on
storage, we prove first-order control and state error estimates uniform
as the Lam\'e parameter increases.
Numerical experiments confirm the predicted convergence rates, examine the
dependence on the Lam\'e and storage coefficients, and demonstrate
localized pressure control with active bounds.
\end{abstract}

\maketitle
\section{Introduction}
\label{sec:intro}

Poroelasticity models the coupled interaction of solid
deformation and pore-fluid flow in a saturated porous medium.
Its mathematical foundations trace back to Biot's theory of
consolidation \cite{biot1941general,wang2000poroelasticity,cheng2016poroelasticity},
which describes the time-dependent
compaction of a porous elastic solid under fluid pressure.
Systems of this kind arise in geomechanics, including subsurface
fluid storage \cite{WANGEN2016486}, in bone mechanics
\cite{Swan200325}, and in consolidation caused by groundwater
withdrawal \cite{Kim1999549}.

Analytical and numerical studies of the linear Biot system include
\cite{Showalter:2000,MR1393902,MR747212}. A central difficulty in its
numerical approximation is \emph{volumetric locking}: when the solid
is nearly incompressible ($\lambda\to\infty$), standard low-order
displacement approximations can lose accuracy. 
Numerical methods to circumvent this difficulty were proposed including stress-velocity mixed formulation \cite{Lee:2016a}, mixed method with reduced integration \cite{yi2017locking}, nonconforming finite element method \cite{MR3803860}, the total pressure formulation
\cite{Oyarzua_Ruiz-Baier:2016,Lee_Mardal_Winther:2017}, a stabilized method \cite{rodrigo2018new}, and different mixed/hybridized-mixed formulations  \cite{hong2019parameter,CLR:2023}, to name a few. This list of locking-free methods is not exhaustive. 

Optimal control of PDE systems is treated, for example, in
\cite{hinze2009optimization}. For poroelasticity, well-posedness of optimal control problems was studied in \cite{MR4410836} with the zero-storage assumption, while numerical methods for distributed-control problem was studied in \cite{khan2026distributed} with general nonnegative storage coefficient.
A recent work \cite{gudi2026virtual} discussed a three-field virtual element method for distributed and Neumann boundary controls for static poroelasticity problems. An extension to thermo-poroelasticity optimal control problems was also studied recently in \cite{Lee:2026a}. 

As observed in \cite{Showalter:2000}, the quasi-static poroelasticity system can be reduced to a parabolic-type evolutionary equation of pressure via mechanical Schur complement. Thus, the analysis of parabolic optimal control problems provides a framework for the space-time discretization and error estimates considered here; see, e.g., \cite{Meidner-Vexler:2008a,Meidner-Vexler:2008b}. For the quasi-static poroelasticity system, however, the regularity of the state and adjoint variables must be established through the saddle point problem coupling the displacement and total pressure. Moreover, eliminating the mechanical variables from the continuous and discrete formulations leads to different effective mass operators. The resulting consistency error cannot be estimated directly under the available time regularity and requires a duality argument combined with discrete summation by parts. The dependence of these estimates on the storage coefficient and the Lamé parameter must also be controlled.

In this work we study the distributed optimal control of a linear poroelastic
system in the total pressure formulation.  The control is a distributed
fluid-pressure source, and the objective functional combines
running-cost tracking of the displacement and fluid pressure,
together with Tikhonov regularization.  The main contributions of this paper are as follows.
\begin{enumerate}
    \item We prove well-posedness of the differentiated total pressure poroelasticity system for nonnegative storage and prove existence and uniqueness of an optimal control. 
    \item By a stability analysis of the system via mechanical Schur complement, we characterize the parameter dependence of the effective mass and show that its uniform coercivity is lost in the simultaneous zero-storage $(s_0 \downarrow 0)$ and incompressible limit ($\lambda \to +\infty$). Sufficient conditions for uniform coercivity of effective mass are also obtained.
    \item We discretize the state and adjoint systems by stable mixed finite elements in space and the dG(0) method in time, and derive first-order a priori control and state error estimates. Adopting the duality argument of parabolic optimal control problems in \cite{Meidner-Vexler:2008a,Meidner-Vexler:2008b} we obtain the error estimates with minimal regularity assumptions of solutions. These error estimates are also robust for parameters under the sufficient conditions for uniform coercivity of effective mass. 
\end{enumerate}

The remainder of the paper is organized as follows.
Section~\ref{sec:wellposedness} introduces the model, derives the
total pressure weak formulation, and establishes well-posedness.
Section~\ref{sec:OCP} formulates the optimal control problem, proves
existence of optimal controls, and derives first-order necessary
optimality conditions.
Section~\ref{sec:THM:FEM} introduces the fully discrete scheme and
the discrete optimality system.
Section~\ref{sec:THM:error} presents the a priori error analysis for the
discrete controls and states.
Section~\ref{sec:numerical} reports convergence, parameter-robustness,
and localized pressure-control experiments.

\section{Problem formulation and well-posedness}
\label{sec:wellposedness}

\subsection{Notation and function spaces}
\label{sec:notations}

Let $\Omega\subset\mathbb{R}^d$, $d\in\{2,3\}$, be a bounded domain
with Lipschitz boundary $\partial\Omega$.  We decompose
$\partial\Omega=\overline{\Gamma_D\cup\Gamma_N}$ with
$\Gamma_D\cap\Gamma_N=\emptyset$, $|\Gamma_D|>0$, and $|\Gamma_N|>0$,
where $\Gamma_D$
is the Dirichlet part of the boundary for the displacement and
$\Gamma_N$ the Neumann (traction-free) part.  The time interval is $(0,T)$ with $0<T<\infty$.

For $k\ge0$ and $1\le r\le\infty$, $W^{k,r}(\Omega)$ denotes the
standard Sobolev space on $\Omega$ with norm $\|\cdot\|_{W^{k,r}}$,
and $H^k(\Omega):=W^{k,2}(\Omega)$ with norm $\|\cdot\|_{H^k}$ and
seminorm $|\cdot|_{H^k}$.  We write $L^2(\Omega):=H^0(\Omega)$ and
denote its inner product by $(\cdot,\cdot)_\Omega$ and norm by
$\|\cdot\|_0:=\|\cdot\|_{L^2(\Omega)}$.  For $k\ge1$ we also abbreviate
$\|\cdot\|_k:=\|\cdot\|_{H^k}$ and $|\cdot|_k:=|\cdot|_{H^k}$.  The space $H^1_0(\Omega)$
is the closure of $C^\infty_c(\Omega)$ in $H^1(\Omega)$.  Its dual
is denoted $H^{-1}(\Omega)$, with duality pairing
$\langle\cdot,\cdot\rangle_{H^{-1},H^1_0}$.  For vector-valued
functions we write $\boldsymbol{H}^k(\Omega):=(H^k(\Omega))^d$ and use
the same notation $(\cdot,\cdot)_\Omega$ for the $L^2$-inner product.
The symmetric gradient is $\varepsilon(v):=\frac{1}{2}(\nabla
v+(\nabla v)^\top)$ for $v\in\boldsymbol{H}^1(\Omega)$.

For a Banach space $X$ and $1\le r\le\infty$, $L^r(0,T;X)$ denotes
the Bochner space of strongly measurable functions $v:(0,T)\to X$
with $\|v\|_{L^r(0,T;X)}<\infty$.  We write $L^2_tL^2_x :=
L^2(0,T;L^2(\Omega))$ for the space-time $L^2$ space, with norm
$\|\cdot\|_{L^2_tL^2_x}^2=\int_0^T\|\cdot\|_0^2\,dt$.  More
generally, $L^r_tX:=L^r(0,T;X)$.  The Bochner--Sobolev space
$W^{k,r}(0,T;X)$ consists of functions whose distributional time
derivatives up to order $k$ belong to $L^r(0,T;X)$; we write
$H^k(0,T;X):=W^{k,2}(0,T;X)$.  The embedding
$H^1(0,T;X)\hookrightarrow C^0([0,T];X)$ holds for any Hilbert space
$X$ \cite{evans1998partial}.  We write $\dot v:=\partial_t v$ for the
time derivative.

\subsection{Governing equations and total pressure formulation}
\label{sec:governing}

In quasi-static poroelasticity equations with homogeneous boundary conditions we seek the solid
displacement $u:(0,T)\times\Omega\to\mathbb{R}^d$ and fluid pressure
$p:(0,T)\times\Omega\to\mathbb{R}$ satisfying 
\begin{subequations}
  \label{eq:tp:strong}
  \begin{align}
    \label{eq:tp:strong1}
    -\operatorname{div}\!\bigl(2\mu\,\varepsilon(u)
      + \lambda\operatorname{div}\,u\,\mathbb{I}
      - \alpha p\,\mathbb{I}\bigr)
    &= \tilde f,
    \\
    \label{eq:tp:strong2}
    s_0\dot p
      + \alpha \operatorname{div}\,\dot u
      - \operatorname{div}(\kappa_p\nabla p)
    &= g_p,
  \end{align}
\end{subequations}
in $(0,T)\times\Omega$, with boundary conditions $u=0$ on
$(0,T)\times\Gamma_D$, zero traction on $(0,T)\times\Gamma_N$, and
$p=0$ on $(0,T)\times\partial\Omega$, together with initial
conditions $(u,p)|_{t=0}=(u_0,p_0)$.  Here
$\mu>0$ and $\lambda>0$ are the Lam\'e parameters,
$\alpha>0$ is the Biot--Willis coupling coefficient,
$s_0\ge0$ is the constrained specific storage coefficient,
and $\kappa_p>0$ is the fluid conductivity. Here $g_p$ 
is the fluid source and $\tilde f$ is the body force. In the distributed optimal control problems that we consider later, $g_p$ will be replaced by the control variable.

We define the function spaces
\begin{align}
  \label{eq:tp:spaces}
  V &:= \bigl\{v\in\boldsymbol{H}^1(\Omega) : v|_{\Gamma_D}=0\bigr\},
  &
  Q_T &:= L^2(\Omega),
  &
  Q &:= H^1_0(\Omega),
\end{align}
equipped with inner products
$(u,v)_V := (2\mu\,\varepsilon(u),\varepsilon(v))_\Omega$,
$(p,q)_{Q_T}:=(p,q)_\Omega$, and
$(p,q)_Q := (\kappa_p\nabla p,\nabla q)_\Omega$,
and induced norms
\begin{align}
  \label{eq:tp:norms}
  \|v\|_V := \sqrt{(v,v)_V},
  \qquad
  \|p\|_Q := \sqrt{(p,p)_Q}
            = (\kappa_p\nabla p,\nabla p)_\Omega^{1/2}.
\end{align}
By Korn's inequality, $\|\cdot\|_V$ is equivalent to
$\|\cdot\|_{\boldsymbol{H}^1}$ on $V$.  By the Poincar\'e inequality,
$\|\cdot\|_Q$ is equivalent to $\|\cdot\|_{H^1}$ on $H^1_0(\Omega)$,
with constant depending on $\kappa_p$.  We denote by $V^*$ the dual of
$V$ and keep the notation $(g,v)_\Omega$ for the $V^*$--$V$ duality
pairing when $g\in V^*$; it coincides with the $L^2$ inner product when
$g\in L^2(\Omega)^d$.  The body force satisfies
$\tilde f\in H^1(0,T;V^*)$.

When $\lambda\gg1$ (nearly incompressible solid), standard
displacement-pressure finite elements applied directly to
\eqref{eq:tp:strong} suffer from \emph{volumetric locking}: the
discrete solution locks onto the constraint $\operatorname{div} u\approx0$
and the approximation deteriorates. To circumvent this, we introduce
the \emph{total pressure}
\begin{align}
  \label{eq:total-pressure-def}
  p_T := -\lambda \operatorname{div} u + \alpha p,
\end{align}
equivalently $\operatorname{div}u = \lambda^{-1}(-p_T+\alpha p)$,
as an additional unknown, as in \cite{Oyarzua_Ruiz-Baier:2016,Lee_Mardal_Winther:2017}.

We introduce the bilinear forms
\begin{align}
  \label{eq:tp:bilinear-forms}
  a_u(u,v) &:= (2\mu\,\varepsilon(u),\varepsilon(v))_\Omega,
  &
  b(q_T,v) &:= -(q_T,\operatorname{div} v)_\Omega,
\end{align}
for $(u,v)\in V\times V$ and $(q_T,v)\in Q_T\times V$.  The form
$a_u$ is coercive on $V$ by Korn's
inequality, and $b$ is the
divergence coupling form arising in mixed formulations. Here the coefficient $s_1$ is defined as
\begin{align}
  \label{eq:tp:s-coefficients}
  s_1 := s_0+\lambda^{-1}\alpha^2.
\end{align}

The storage coupling of the total pressure and fluid pressure is represented by
\begin{align}
  \label{eq:tp:S-matrix}
  \mathbb{S}(\lambda) :=
  \begin{bmatrix}
    \lambda^{-1} & -\lambda^{-1}\alpha \\
    -\lambda^{-1}\alpha & s_1
  \end{bmatrix}.
\end{align}
The associated quadratic form is precisely the storage contribution
that appears when the total pressure constraint and fluid mass balance
are combined. Its nonnegativity, established by the following lemma,
is useful for identifying the parameter dependence of the storage block.

\begin{lemma}[$\lambda$-uniform nonnegativity of $\mathbb{S}$]
  \label{lem:tp:S-nonneg}
  For $\lambda>0$ and $s_0\ge0$,
  \begin{align}
    \label{eq:tp:S-decomposition}
    \mathbb S(\lambda)=\lambda^{-1}ww^\top
       +\operatorname{diag}(0,s_0),\qquad
       w=(1,-\alpha)^\top,
  \end{align}
  and hence, for all $\xi=(\xi_{p_T},\xi_p)^\top$,
  \begin{align}
    \label{eq:tp:S-quadratic-form}
    \xi^\top\mathbb S(\lambda)\xi
      &=\lambda^{-1}(\xi_{p_T}-\alpha\xi_p)^2+s_0\xi_p^2,\\
    \label{eq:tp:S-lower-bound}
    \xi^\top\mathbb S(\lambda)\xi &\ge s_0\xi_p^2.
  \end{align}
  In particular the lower bound is independent of $\lambda>0$.
\end{lemma}

Differentiating the momentum balance \eqref{eq:tp:strong1} in time
(with $f:=\dot{\tilde f}$) and using $\dot p_T=-\lambda\operatorname{div}\dot u
+\alpha\dot p$ to eliminate $\operatorname{div}\dot u$ from
\eqref{eq:tp:strong2} yields the
\emph{differentiated total pressure system}
\begin{subequations}
  \label{eq:tp:strong-pT}
  \begin{align}
    -\operatorname{div}(2\mu\,\varepsilon(\dot u) - \dot p_T\,\mathbb{I})
    &= f,
    \\
    -\operatorname{div}\dot u
      - \lambda^{-1}\dot p_T
      + \lambda^{-1}\alpha\dot p
    &= 0,
    \\
    -\lambda^{-1}\alpha\dot p_T
      + s_1\dot p
      - \operatorname{div}(\kappa_p\nabla p)
    &= g_p,
  \end{align}
\end{subequations}
with $s_1$ defined in \eqref{eq:tp:s-coefficients}.
All factors of $\lambda$ now appear only as $\lambda^{-1}$, so
\eqref{eq:tp:strong-pT} is a first-order evolution problem
for $(u,p_T,p)$ with coefficients bounded uniformly as
$\lambda\to\infty$.  The precise equivalence between
\eqref{eq:tp:strong} and \eqref{eq:tp:strong-pT} under compatible
initial data is established in \S\ref{sec:equivalence}.

The continuous inf-sup condition on the full total pressure space is
\begin{align}
  \label{eq:tp:inf-sup}
  \inf_{0\ne q_T\in Q_T}\sup_{0\ne v\in V}
  \frac{(q_T,\operatorname{div}v)_\Omega}
       {\|v\|_1\|q_T\|_0}\ge\beta>0.
\end{align}
Since $|\Gamma_N|>0$, the divergence operator
$\operatorname{div}:V\to L^2(\Omega)=Q_T$ is surjective.

\begin{rmk}
  \label{rmk:tp:lambda-indep-babuska}
  The forms $a_u$ and $b$, and the inf-sup constant in
  \eqref{eq:tp:inf-sup}, are independent of $\lambda$ and $s_0$.
  The parameter dependence of the evolution estimates is determined
  separately by the effective mass introduced below.
\end{rmk}

The continuous weak formulation of \eqref{eq:tp:strong-pT} reads:
find $(u,p_T,p)\in V\times Q_T\times Q$ such that, for a.e.\ $t\in(0,T)$,
\begin{subequations}
  \label{eq:tp:weak-state}
  \begin{align}
    \label{eq:tp:G1}
    a_u(\dot u,v)+b(\dot p_T,v)&=(f,v)_\Omega,\\
    \label{eq:tp:G2}
    b(q_T,\dot u)-\lambda^{-1}(\dot p_T,q_T)_\Omega
      +\lambda^{-1}\alpha(\dot p,q_T)_\Omega&=0,\\
    \label{eq:tp:G3}
    -\lambda^{-1}\alpha(\dot p_T,q)_\Omega+s_1(\dot p,q)_\Omega
      +(\kappa_p\nabla p,\nabla q)_\Omega&=(g_p,q)_\Omega
  \end{align}
\end{subequations}
for all $(v,q_T,q)\in V\times Q_T\times Q$.
The time regularity and well-posedness of this system are established
in Sections~\ref{sec:apriori}--\ref{sec:existence}.

\subsection{Equivalence with the original formulation}
\label{sec:equivalence}

The weak system \eqref{eq:tp:weak-state} is obtained by differentiating
momentum balance and the total pressure relation in time. We establish
its equivalence with the original system in the continuous spaces
$V,Q_T=L^2(\Omega),Q$.

Define, for $(u,p_T,p) \in C^1(0,T;V\times L^2(\Omega)\times Q)$, the \emph{total pressure residual}
\begin{align}
  \label{eq:tp:Z-residual}
  Z(t) := p_T(t) + \lambda \operatorname{div}u(t)
    - \alpha p(t) \in L^2(\Omega),
\end{align}
so that $Z(t)=0$ is precisely the algebraic total pressure relation
\eqref{eq:total-pressure-def}. 
The differentiated formulation permits arbitrary initial data;
compatibility is imposed when equivalence with the original Biot system
is asserted.

The \emph{original system}
consists of the (non-differentiated) momentum and mass-balance
equations: find $(u,p)\in C^1(0,T;V\times Q)$ such that
\begin{subequations}
  \label{eq:tp:original}
  \begin{align}
    \label{eq:tp:O1}
    a_u(u,v) + \lambda(\operatorname{div}u,\operatorname{div}v)_\Omega
    - \alpha(p,\operatorname{div} v)_\Omega
    &= (\tilde f, v)_\Omega,
    \\
    \label{eq:tp:O2}
    s_0(\dot p,q)_\Omega
    + \alpha(\operatorname{div}\dot u,q)_\Omega
    + (\kappa_p\nabla p,\nabla q)_\Omega
    &= (g_p,q)_\Omega,
  \end{align}
\end{subequations}
for all $(v,q)\in V\times Q$ and all $0<t<T$,
where $\tilde f \in H^1(0,T;V^*)$ and $f:=\dot{\tilde f}$ is the
forcing appearing in \eqref{eq:tp:G1}. We say the initial data
$(u(0),p_T(0),p(0))$ are \emph{compatible} if both of the following hold:
\begin{subequations}
\label{eq:tp:compatible-data}
\begin{align}
  \label{eq:tp:compatible-data-pT}
  Z(0) = p_T(0) + \lambda \operatorname{div}u(0)
    - \alpha p(0) = 0,
\end{align}
i.e., the algebraic total pressure relation \eqref{eq:total-pressure-def} at $t=0$, and
\begin{align}
  \label{eq:tp:compatible-data-mech}
  a_u(u(0),v) + b(p_T(0),v) = (\tilde f(0),v)_\Omega
  \qquad \forall\, v\in V,
\end{align}
\end{subequations}
i.e., the momentum equation holds at $t=0$ with the initial total pressure
$p_T(0)=-\lambda\operatorname{div}u(0)+\alpha p(0)$.
Condition \eqref{eq:tp:compatible-data-mech} is needed to recover the
original (non-differentiated) momentum equation from its time-integrated
form. It is automatically satisfied if $(u(0),p(0))$ solves the
original static equilibrium at $t=0$.

\begin{proposition}[Equivalence of the original and differentiated formulations]
  \label{prop:tp:equivalence}
  Let $V,Q_T=L^2(\Omega),Q=H^1_0(\Omega)$ be the
  continuous function spaces of \eqref{eq:tp:spaces}.  Assume
  $\tilde f\in H^1(0,T;V^*)$ and $f=\dot{\tilde f}$.
  \begin{enumerate}
  \item[\textup{(i)}] \textup{(Differentiated $\Rightarrow$
    original.)} Let $(u,p_T,p)\in C^0([0,T];V\times
    L^2(\Omega)\times Q)\cap C^1(0,T;V\times
    L^2(\Omega)\times Q)$ solve \eqref{eq:tp:weak-state}
    with compatible initial data \eqref{eq:tp:compatible-data}. Then
    $Z(t)=0$ for all $t\in[0,T]$, and $(u,p)$ solves the
    original system \eqref{eq:tp:original} with
    \begin{align}
      \label{eq:tp:ftilde-reconstruction}
      \tilde f(t) := \tilde f(0) + \int_0^t f(s)\,ds.
    \end{align}

  \item[\textup{(ii)}] \textup{(Original $\Rightarrow$
    differentiated.)} Conversely, let $(u,p) \in
    C^2(0,T;V\times Q)$ solve \eqref{eq:tp:original}, and
    define $p_T(t):=-\lambda\,\operatorname{div}\,u(t)+\alpha p(t)$
    via \eqref{eq:total-pressure-def} (so $Z(t)\equiv0$). Then
    $(u,p_T,p)\in C^1(0,T;V\times L^2(\Omega)\times Q)$
    solves \eqref{eq:tp:weak-state} with $f=\dot{\tilde f}$.
  \end{enumerate}
\end{proposition}

\begin{proof}
  \textbf{(i)} Test the $p_T$-equation of \eqref{eq:tp:weak-state}
  with an arbitrary $q_T\in L^2(\Omega)$
  (this is valid at the continuous level where $Q_T=L^2(\Omega)$):
  \begin{align*}
    -(\operatorname{div}\,\dot u,q_T)_\Omega
    - \lambda^{-1}(\dot p_T,q_T)_\Omega
    + \lambda^{-1}\alpha(\dot p,q_T)_\Omega = 0
    \qquad \forall\,q_T\in L^2(\Omega).
  \end{align*}
  Since this holds for all $q_T\in L^2(\Omega)$, it holds pointwise
  a.e.\ in $\Omega$:
  \begin{align*}
    \dot p_T = -\lambda \operatorname{div}\dot u + \alpha\dot p
    = \frac{d}{dt}\bigl(-\lambda \operatorname{div}u+\alpha p\bigr).
  \end{align*}
  Hence $\dot Z(t)=0$ for all $t\in(0,T)$, so $Z(t)=Z(0)=0$ by
  \eqref{eq:tp:compatible-data}, i.e., $p_T(t)=-\lambda\operatorname{div}u(t)+\alpha p(t)$
  for all $t$.  Substituting $\dot p_T=-\lambda\operatorname{div}\dot u+\alpha\dot p$
  into the displacement equation \eqref{eq:tp:G1} and integrating in time from $0$ to $t$ gives
  \begin{align*}
    a_u(u(t),v) + b(p_T(t),v)
    = a_u(u(0),v) + b(p_T(0),v)
    + \int_0^t (f(s),v)_\Omega\,ds.
  \end{align*}
  By the initial momentum condition \eqref{eq:tp:compatible-data-mech},
  $a_u(u(0),v)+b(p_T(0),v)=(\tilde f(0),v)_\Omega$, so
  \begin{align*}
    a_u(u(t),v)+b(p_T(t),v)
    = \Bigl(\tilde f(0)+\int_0^t f(s)\,ds,v\Bigr)_\Omega
    = (\tilde f(t),v)_\Omega,
  \end{align*}
  which recovers \eqref{eq:tp:O1}.  Substituting into \eqref{eq:tp:G3} and using $s_1=s_0+\lambda^{-1}\alpha^2$
  (cf.\ \eqref{eq:tp:s-coefficients}) recovers \eqref{eq:tp:O2}.

  \textbf{(ii)} Conversely, with $Z\equiv0$ by definition of $p_T$,
  differentiating \eqref{eq:tp:O1} in time gives \eqref{eq:tp:G1}. Substituting
  $\operatorname{div}\dot u=\lambda^{-1}(-\dot p_T+\alpha\dot p)$
  into \eqref{eq:tp:O2} and using the $s$-coefficient
  identity gives 
  \eqref{eq:tp:G3}; and $Z\equiv0$ implies 
  \eqref{eq:tp:G2}. Hence $(u,p_T,p)$ solves
  \eqref{eq:tp:weak-state}.
\end{proof}

\begin{rmk}[Compatible initial data]
  \label{rmk:tp:compatibility}
  The differentiated total pressure system is equivalent to the
  original poroelastic system under both parts of
  \eqref{eq:tp:compatible-data}: the algebraic condition
  \eqref{eq:tp:compatible-data-pT} ensures $Z(t)\equiv0$, while
  the initial momentum condition \eqref{eq:tp:compatible-data-mech}
  is needed to identify the integration constant when recovering the
  original momentum equation~\eqref{eq:tp:O1}.
  Under \eqref{eq:tp:compatible-data-pT}, $p_T(0)$ is determined by $u(0)$ and $p(0)$
  through \eqref{eq:total-pressure-def}, and then
  \eqref{eq:tp:compatible-data-mech} is automatically satisfied
  if $(u(0),p(0))$ solves the original static momentum balance at $t=0$.
\end{rmk}

\subsection{A priori estimates and well-posedness via mechanical Schur complement}
\label{sec:apriori}

For a priori estimate we eliminate
the differentiated mechanical variables $(\dot u,\dot p_T)$ and use the
elastic Schur complement as an additional contribution to the scalar
mass form. This allows an estimate independent of $s_0^{-1}$,  

Since the divergence map $\operatorname{div}:V\to
L^2(\Omega)$ is surjective, there exists a bounded right inverse
$R:L^2(\Omega)\to V$ such that
\begin{align}
    \label{eq:R-operator}
  \operatorname{div}(Rq)=q,
  \qquad \|Rq\|_V\le C_D\|q\|_0 .
\end{align}

\begin{lemma}[auxiliary solution operator]
  \label{lem:tp:mech-operator}
  For given $p\in L^2(\Omega)$, let
  $(\mathcal{U}p,\mathcal{P}_Tp)\in V\times Q_T$ be the unique
  solution of
  \begin{subequations}      
    \label{eq:mech-op-cont}
  \begin{align}
    \label{eq:tp:mech-op-cont1}
    a_u(\mathcal{U}p,v) + b(\mathcal{P}_Tp,v) &= 0,
    \\
    \label{eq:tp:mech-op-cont2}
    b(q_T,\mathcal{U}p) - \lambda^{-1}(\mathcal{P}_Tp,q_T)
    + \lambda^{-1}\alpha(p,q_T)
    &= 0
  \end{align}
  \end{subequations}
  for all $(v,q_T)\in V\times Q_T$.  Then
  \begin{align}
    \label{eq:tp:mech-op-bound}
    \|\mathcal{U}p\|_1 + \|\mathcal{P}_Tp\|_0
    \le C\|p\|_0,
  \end{align}
  where $C$ is independent of $s_0$ and can be chosen uniformly for
  $\lambda$ in any interval $[\lambda_0,\infty)$, $\lambda_0>0$.
  For $p\in H^1_0(\Omega)$, additionally
  $\|\mathcal{P}_Tp\|_{H^1}\le C\|p\|_{H^1}$, and for
  $p\in H^2(\Omega)$,
  $\|\mathcal{P}_Tp\|_{H^1}\le C\|p\|_{H^2}$, under the regularity
  assumption~(\ref{item:app:stokes}) of Appendix~\ref{app:regularity},
  applied with $F=0$ and $\zeta=-\lambda^{-1}\alpha p$.
\end{lemma}

\begin{proof}
  This is the saddle point problem of displacement and total pressure with
  right-hand side $-\lambda^{-1}\alpha(p,\cdot)$ in the second equation.
  The Babu\v{s}ka--Brezzi theorem and the inf-sup condition
  give existence, uniqueness, and \eqref{eq:tp:mech-op-bound}.  
\end{proof}

\begin{lemma}
  \label{lem:tp:reduced-mass}
  For $\mathcal{U}$ and $\mathcal{P}_T$ in \eqref{eq:mech-op-cont}, define
  \begin{align}
    \label{eq:tp:reduced-mass-def}
    m_{\rm eff}(p,q):= -\lambda^{-1}\alpha(\mathcal{P}_Tp,q)_{\Omega}
      + s_1(p,q)_{\Omega}.
  \end{align}
  Then $m_{\rm eff}(\cdot,\cdot)$ is symmetric and bounded.  Moreover, for every $s_0\ge0$,
  \begin{align}
    \label{eq:tp:reduced-mass-identity}
    m_{\rm eff}(q,q)
    &= s_0\|q\|_0^2
      + a_u(\mathcal Uq,\mathcal Uq)
      + \lambda^{-1}\|\mathcal P_Tq-\alpha q\|_0^2
     \ge c_{\rm eff} \|q\|_0^2    
  \end{align}
  where $c_{\rm eff}:= s_0 + \frac{\alpha^2}{\lambda + C_D^2}$.
\end{lemma}

\begin{proof}
  Symmetry of $m_{\rm eff}(\cdot,\cdot)$ follows by testing the two mixed problems \eqref{eq:mech-op-cont} for
  $(\mathcal Up,\mathcal P_Tp)$ and
  $(\mathcal Uq,\mathcal P_Tq)$ against one another. Testing \eqref{eq:tp:mech-op-cont1} with $\mathcal Up$ and \eqref{eq:tp:mech-op-cont2} with $\mathcal P_Tp$ gives
  \begin{align*}
    a_u(\mathcal Up,\mathcal Up)
    =\lambda^{-1}\bigl(\alpha(p,\mathcal P_Tp)_{\Omega}
                       -\|\mathcal P_Tp\|_0^2\bigr).
  \end{align*}
  Expanding the square in \eqref{eq:tp:reduced-mass-identity} proves
  the identity in \eqref{eq:tp:reduced-mass-identity}. 
  
  For the inequality in \eqref{eq:tp:reduced-mass-identity}, note that \eqref{eq:tp:mech-op-cont2} also gives   $(\mathcal P_Tp-\alpha p,q_T)_{\Omega}=-\lambda(\operatorname{div}(\mathcal Up), q_T)_{\Omega}$ for every $q_T \in L^2(\Omega)$. Combining this with $q_T=\operatorname{div}v$ and
  \eqref{eq:tp:mech-op-cont1} gives
  \[
    a_u(\mathcal Up,v)
    +\lambda(\operatorname{div}\mathcal Up,\operatorname{div}v)_\Omega
    =\alpha(p,\operatorname{div}v)_\Omega
    \qquad\forall v\in V.
  \]
  Choose $v=Rp$.  By Cauchy--Schwarz in the elastic energy and the
  boundedness of $R$,
  \begin{align*}
    \alpha\|p\|_0^2
    &\le
    \bigl(a_u(\mathcal Up,\mathcal Up)
      +\lambda\|\operatorname{div}\mathcal Up\|_0^2\bigr)^{1/2} \bigl(\|Rp\|_V^2 + \lambda \|p\|_0^2\bigr)^{1/2}
      \\
      &\le \bigl(a_u(\mathcal Up,\mathcal Up)
      +\lambda\|\operatorname{div}\mathcal Up\|_0^2\bigr)^{1/2} \bigl(C_D^2 + \lambda \bigr)^{1/2} \|p\|_0 .
  \end{align*}
  Hence
  $a_u(\mathcal Up,\mathcal Up)
      +\lambda\|\operatorname{div}\mathcal Up\|_0^2
    \ge \frac{\alpha^2}{\lambda + C_D^2} \|p\|_0^2$,
  and this implies the last inequality in
  \eqref{eq:tp:reduced-mass-identity}.
\end{proof}

For the differentiated body force $f(t)$, let
$(u_f(t),p_{T,f}(t))\in V\times Q_T$ be the unique solution of
\begin{align}
  \label{eq:tp:mechanical-load-lift}
  a_u(u_f,v)+b(p_{T,f},v)&=(f,v)_\Omega,
  \notag\\
  b(q_T,u_f)-\lambda^{-1}(p_{T,f},q_T)_{\Omega}&=0
\end{align}
for all $(v,q_T)\in V\times Q_T$.  The same saddle-point argument gives
\begin{align}
  \label{eq:tp:mechanical-load-lift-bound}
  \|u_f\|_V+\|p_{T,f}\|_0\le C\|f\|_0,
  \qquad p_{T,f}=-\lambda\operatorname{div}u_f.
\end{align}

\begin{lemma}
  \label{lem:tp:schur-reduction}
  A sufficiently regular triple $(u,p_T,p)$ solves
  \eqref{eq:tp:G1}--\eqref{eq:tp:G3} if and only if
  \begin{align}
    \label{eq:tp:mechanical-reconstruction}
    \dot u=u_f+\mathcal U\dot p,
    \qquad
    \dot p_T=p_{T,f}+\mathcal P_T\dot p,
  \end{align}
  and $p$ solves the reduced scalar problem
  \begin{align}
    \label{eq:tp:reduced-continuous}
    m_{\rm eff}(\dot p,q)+(\kappa_p\nabla p,\nabla q)_\Omega
    =(g_p,q)_\Omega+\lambda^{-1}\alpha(p_{T,f},q)_\Omega
    \qquad\forall q\in Q.
  \end{align}
\end{lemma}

\begin{proof}
  Equations \eqref{eq:tp:G1}--\eqref{eq:tp:G2} and linearity of the
  mixed mechanical problem give \eqref{eq:tp:mechanical-reconstruction}.
  Substitution in \eqref{eq:tp:G3} gives
  \eqref{eq:tp:reduced-continuous}.  Conversely, reconstructing
  $(\dot u,\dot p_T)$ by \eqref{eq:tp:mechanical-reconstruction}
  gives \eqref{eq:tp:G1}--\eqref{eq:tp:G2}, and
  \eqref{eq:tp:reduced-continuous} is exactly \eqref{eq:tp:G3}.
\end{proof}

\begin{lemma}
  \label{lem:tp:energy}
  Let $s_0\ge0$, and let $(u,p_T,p)$ solve the differentiated
  total pressure system with $p(0)\in Q$.  For $0<t\le T$,
  \begin{align}
    \label{eq:tp:energy2}
    &\|\dot p\|_{L^2(0,t;L_x^2)}
      +\|p\|_{C^0([0,t];Q)}
      +\|\dot u\|_{L^2(0,t;V)}
      +\|\dot p_T\|_{L^2(0,t;L_x^2)}
    \notag\\
    &\qquad\le C_{\rm eff}\bigl(
      \|p(0)\|_Q+\|f\|_{L^2(0,t;L_x^2)}
      +\|g_p\|_{L^2(0,t;L_x^2)}\bigr),\\
    \label{eq:tp:energy1}
    &\|p\|_{C^0([0,t];L_x^2)}+\|p\|_{L^2(0,t;Q)}
      +\|u\|_{C^0([0,t];V)}+\|p_T\|_{C^0([0,t];L_x^2)}
    \notag\\
    &\qquad\le \|u(0)\|_V+\|p_T(0)\|_0+C\|p(0)\|_0
      +C\|f\|_{L^1(0,t;L_x^2)}
    \notag\\
    &\qquad\quad+C_{\rm eff}\bigl(
      m_{\rm eff}(p(0),p(0))^{1/2}
      +\|f\|_{L^2(0,t;L_x^2)}
      +\|g_p\|_{L^2(0,t;L_x^2)}\bigr).
  \end{align}
  Here $C_{\rm eff}$ depends on $c_{\rm eff}^{-1}$ in
  \eqref{eq:tp:reduced-mass-identity}. 
\end{lemma}

\begin{proof}
  Test \eqref{eq:tp:reduced-continuous} first with $q=\dot p$.
  Since $m_{\rm eff}$ is coercive and
  $(\kappa_p\nabla p,\nabla\dot p)_\Omega=\tfrac12\frac{d}{dt}\|p\|_Q^2$,
  the mechanical lift estimate and Young's inequality yield
  \eqref{eq:tp:energy2}, using
  \eqref{eq:tp:mechanical-reconstruction} for $\dot u$ and $\dot p_T$.

  If we test \eqref{eq:tp:reduced-continuous} with $q=p$, then symmetry of $m_{\rm eff}$ gives
  \[
    \frac12\frac{d}{dt}m_{\rm eff}(p,p)+(\kappa_p\nabla p,\nabla p)_\Omega
      =(g_p,p)_\Omega+\lambda^{-1}\alpha(p_{T,f},p)_\Omega.
  \]
  Poincar\'e's inequality and
  $\lambda^{-1}p_{T,f}=-\operatorname{div}u_f$ give
  \[
    \frac{d}{dt}m_{\rm eff}(p,p)+c\|p\|_Q^2
    \le C\bigl(\|g_p\|_0^2+\|f\|_0^2\bigr).
  \]
  Integrating and using $m_{\rm eff}(p,p)\ge c_{\rm eff}\|p\|_0^2$ yields
  \[
    \|p\|_{C^0([0,t];L_x^2)}+\|p\|_{L^2(0,t;Q)}
    \le C_{\rm eff}\bigl(m_{\rm eff}(p(0),p(0))^{1/2}
       +\|f\|_{L^2(0,t;L_x^2)}+\|g_p\|_{L^2(0,t;L_x^2)}\bigr).
  \]
  Finally integrate the mechanical reconstruction in time:
  \[
    u(s)=u(0)+\int_0^s u_f(r)\,dr+\mathcal U(p(s)-p(0)),
    \qquad
    p_T(s)=p_T(0)+\int_0^s p_{T,f}(r)\,dr
       +\mathcal P_T(p(s)-p(0)).
  \]
  The mechanical operator and load-lift bounds yield
  \[
    \|u\|_{C^0([0,t];V)}+\|p_T\|_{C^0([0,t];L_x^2)}
    \le\|u(0)\|_V+\|p_T(0)\|_0
       +C\|f\|_{L^1(0,t;L_x^2)}
       +C(\|p\|_{C^0([0,t];L_x^2)}+\|p(0)\|_0),
  \]
  and thus \eqref{eq:tp:energy1}.
\end{proof}

\subsection{Existence and uniqueness of weak solutions}
\label{sec:existence}

\begin{definition}
  \label{def:tp:weak-solution}
  For given $(u_0,p_{T,0},p_0)\in V\times L^2(\Omega)\times Q$,
  a triple
  \[
    (u,p_T,p)\in H^1(0,T;V)\times H^1(0,T;L^2(\Omega))
      \times\bigl(H^1(0,T;L^2(\Omega))\cap L^2(0,T;Q)\bigr)
  \]
  is a weak solution if
  \begin{subequations}
    \label{eq:tp:weak-solution}
    \begin{align}
      \label{eq:tp:ws1}
      \int_0^T a_u(\dot u,v)\,dt+\int_0^T b(\dot p_T,v)\,dt
      &=\int_0^T(f,v)_\Omega\,dt,\\
      \label{eq:tp:ws2}
      \int_0^T\bigl[b(q_T,\dot u)-\lambda^{-1}(\dot p_T,q_T)_\Omega
      +\lambda^{-1}\alpha(\dot p,q_T)_\Omega\bigr]dt&=0,\\
      \label{eq:tp:ws3}
      \int_0^T\bigl[-\lambda^{-1}\alpha(\dot p_T,q)_\Omega
      +s_1(\dot p,q)_\Omega+(\kappa_p\nabla p,\nabla q)_\Omega\bigr]dt
      &=\int_0^T(g_p,q)_\Omega\,dt,
    \end{align}
  \end{subequations}
  for all temporally compactly supported test functions
  $(v,q_T,q) \in C_0^{\infty}\bigl((0,T); V \times Q_T \times Q \bigr)$.
\end{definition}

\begin{theorem}
  \label{thm:tp:existence}
  Let $s_0\ge0$, $f,g_p\in L^2(0,T;L^2(\Omega))$, and
  $(u_0,p_{T,0},p_0)\in V\times L^2(\Omega)\times Q$.  Then there exists a
  unique weak solution in the sense of Definition~\ref{def:tp:weak-solution},
  and it satisfies \eqref{eq:tp:energy2}--\eqref{eq:tp:energy1}.
\end{theorem}

\begin{proof}
  By Lemma~\ref{lem:tp:reduced-mass}, $m_{\rm eff}$ is a bounded symmetric
  coercive bilinear form on $L^2(\Omega)$ for every fixed
  $\lambda>0$ and every $s_0\ge0$.  The form $(\kappa_p\nabla\cdot,\nabla\cdot)_\Omega$ is coercive
  on $Q=H^1_0(\Omega)$.  Therefore the standard theory for linear
  parabolic problems applied to \eqref{eq:tp:reduced-continuous} gives
  a unique
  \[
    p\in H^1(0,T;L^2(\Omega))\cap L^\infty(0,T;Q)
  \]
  with $p(0)=p_0$.  Define $\dot u$ and $\dot p_T$ by
  \eqref{eq:tp:mechanical-reconstruction} and set
  \[
    u(t)=u_0+\int_0^t\dot u(s)\,ds,
    \qquad
    p_T(t)=p_{T,0}+\int_0^t\dot p_T(s)\,ds.
  \]
  Lemma~\ref{lem:tp:mech-operator} and
  \eqref{eq:tp:mechanical-load-lift-bound} imply
  $u\in H^1(0,T;V)$ and $p_T\in H^1(0,T;L^2(\Omega))$.
  Lemma~\ref{lem:tp:schur-reduction} then shows that the reconstructed
  triple satisfies the three-field differentiated system and hence the
  weak formulation.  Uniqueness follows first for $p$ from the reduced
  parabolic problem and then for $(\dot u,\dot p_T)$ from the unique
  mechanical reconstruction; the prescribed initial data determine
  $(u,p_T)$ uniquely.
\end{proof}

\section{Optimal control problem}
\label{sec:OCP}

\subsection{Existence of optimal control}
\label{sec:OCP:existence}

Let $l_p < L_p$ be prescribed lower and upper control bounds.  The
\emph{admissible control set} is
\begin{align}
  \label{eq:tp:admissible-set}
  \mathcal{M}_{p,ad}
  := \bigl\{
    m \in L^2(0,T;L^2(\Omega)) :
    l_p \le m(t,x) \le L_p
    \ \text{a.e.\ in } (0,T)\times\Omega
  \bigr\}.
\end{align}
The set $\mathcal{M}_{p,ad}$ is nonempty (the constant $l_p$ belongs to it),
closed and convex in $L^2(0,T;L^2(\Omega))$, and bounded.

Our main goal is to choose a distributed control
$m_p \in \mathcal{M}_{p,ad}$
acting on the fluid pressure equation such that
the corresponding state $(u, p_T, p)$ best approximates
prescribed target fields $u_C$, $p_C$. 
Let $\alpha_{u, C}$, $\alpha_{p, C} \ge 0$
be tracking weights (with $\alpha_{u,C}+\alpha_{p,C}>0$)
and let $\gamma_p > 0$ be the control cost parameter.
For given target fields
$u_C, p_C \in L_t^2 L_x^2$
and an admissible tuple $\LRp{u, p_T, p, m_p}$,
we define the objective functional
$J \LRp{u, p, m_p}$ as
\begin{align}
	\label{obj_fun}
	J \LRp{u, p, m_p}
	= 
	\frac{\alpha_{u, C}}{2} \| u - u_C\|_{L_t^2L_x^2}^2
	+ \frac{\alpha_{p, C}}{2} \| p - p_C\|_{L_t^2L_x^2}^2
	+ 
	\frac{\gamma_p}{2} \| m_p\|_{L_t^2L_x^2}^2.
\end{align}
The running tracking terms penalize the distance of the trajectory
$(u,p)$ to the prescribed target profiles $(u_C,p_C)$
over the whole time interval $(0,T)$, and the last term is the
cost of the control $m_p$.

\begin{definition}
  \label{def:tp:optimal-control}
  Let $\widehat{m}_p \in \mathcal{M}_{p,ad}$.
  The admissible tuple
  $(\widehat{u}, \widehat{p}_T, \widehat{p}, \widehat{m}_p)$
  is called an \emph{optimal control} if
  \begin{align}
    \label{eq:tp:OCP-optimality}
    J(u, p, m_p)
    \ge
    J(\widehat{u}, \widehat{p}, \widehat{m}_p)
  \end{align}
  for all admissible
  $(u, p_T, p, m_p)$
  with the same initial data
  $(u_0, p_{T,0}, p_0)$.
\end{definition}

Writing $(u(m_p),p(m_p))$ for the displacement and fluid pressure of
the state associated with $m_p$ (for fixed initial data and body force),
we use the reduced functional
\begin{align}
  \label{eq:tp:reduced-functional}
  j(m_p):=J(u(m_p),p(m_p),m_p),\qquad m_p\in\mathcal M_{p,ad}.
\end{align}

\begin{theorem}
  \label{thm:tp:OCP-existence}
  For any given initial data
  $(u_0, p_{T,0}, p_0) \in V \times Q_T \times Q$
  there exists a unique optimal control
  $\widehat{m}_p \in \mathcal{M}_{p,ad}$.
\end{theorem}

\begin{proof}
  By Theorem~\ref{thm:tp:existence}, the solution of the
  differentiated state system depends affinely and continuously on
  $m_p\in L^2(0,T;L^2(\Omega))$, for fixed initial data and body force.
  Denote its displacement and fluid-pressure components by
  $(u(m_p),p(m_p))$ and let $j$ be the reduced functional
  \eqref{eq:tp:reduced-functional}.
  The admissible set $\mathcal M_{p,ad}$ is nonempty, closed and convex,
  hence weakly closed in $L^2(0,T;L^2(\Omega))$.

  The control regularization with $\gamma_p>0$ makes $j$ coercive;
  thus a minimizing sequence is bounded in $L^2(0,T;L^2(\Omega))$
  and has a weakly convergent subsequence with admissible limit.
  The bounded affine state map preserves weak convergence.
  The tracking terms and the control cost are squared Hilbert-space
  norms of affine functions, and are weakly lower semicontinuous.
  Consequently the weak limit minimizes $j$.
  Finally, the control cost is strictly convex because $\gamma_p>0$,
  while the tracking terms are convex.  Hence $j$ is strictly convex
  and its minimizer is unique.
\end{proof}

\subsection{Necessary optimality conditions}
\label{sec:OCP:optimality}

We define the state residual as follows.  For
$(u,p_T,p)$, a control $m_p$, and a test triple
$(v^*,q_T^*,q^*)\in V\times Q_T\times Q$, define
\begin{align}
  \label{eq:tp:residual}
  &\mathcal{R}\bigl((u,p_T,p),m_p;(v^*,q_T^*,q^*)\bigr)
  \notag\\
  &:={}
    a_u(\dot u,v^*) + b(\dot p_T,v^*) + b(q_T^*,\dot u)
    - \lambda^{-1}(\dot p_T,q_T^*)_\Omega
    + \lambda^{-1}\alpha(\dot p,q_T^*)_\Omega
  \notag\\
  &\qquad
    -\lambda^{-1}\alpha(\dot p_T,q^*)_\Omega
    + s_1(\dot p,q^*)_\Omega
    + (\kappa_p\nabla p,\nabla q^*)_\Omega
    - (f,v^*)_\Omega - (m_p,q^*)_\Omega.
\end{align}
Thus the weak state equations \eqref{eq:tp:G1}--\eqref{eq:tp:G3} are
equivalent to
$\mathcal{R}((u,p_T,p),m_p;(v^*,q_T^*,q^*))=0$ for every
$(v^*,q_T^*,q^*)\in V\times Q_T\times Q$.

Using the residual $\mathcal{R}$ of \eqref{eq:tp:residual}, we
introduce the Lagrangian functional
\begin{align*}
  \mathcal{L}\bigl(
    (u, p_T, p),
    (v^*, q_T^*, q^*),
    m_p
  \bigr)
  := J(u, p, m_p)
  - \int_0^T \mathcal{R}\bigl((u,p_T,p),m_p;
      (v^*,q_T^*,q^*)\bigr)\,dt,
\end{align*}
where $(v^*, q_T^*, q^*)$ are the adjoint (co-state) variables
corresponding to the displacement, total pressure, and fluid pressure,
respectively. 

Setting the variation of $\mathcal{L}$ with respect to the state
$(u,p_T,p)$ to zero and integrating the time-derivative terms by parts
yields the adjoint system, whose
time-derivative block is the algebraic transpose of that in the
state equation:
find $(v^*, q_T^*, q^*) \in V \times Q_T \times Q$
such that
\begin{subequations}
\label{eq:tp:adjoint}
\begin{align}
  \label{eq:tp:adj1}
  -a_u(\dot{v}^*, v) - b(\dot{q}_T^*, v)
  &= \alpha_{u,C}(u - u_C, v)_\Omega,
  \\
  \label{eq:tp:adj2}
  -b(q_T, \dot{v}^*) + \lambda^{-1}(\dot{q}_T^*, q_T)_\Omega
  + \lambda^{-1}\alpha (\dot{q}^*, q_T)_\Omega
  &= 0,
  \\
  \label{eq:tp:adj3}
  -\lambda^{-1}\alpha (\dot{q}_T^*, q)_\Omega
  - s_1(\dot{q}^*, q)_\Omega
  + (\kappa_p \nabla q^*, \nabla q)_\Omega
  &= \alpha_{p,C}(p - p_C, q)_\Omega,
\end{align}
\end{subequations}
for all $(v, q_T, q) \in V \times Q_T \times Q$
and for all $0 < t < T$, with the homogeneous terminal condition
\begin{align}
  \label{eq:tp:adj:TC}
  (v^*(T), q_T^*(T), q^*(T)) = (0,0,0).
\end{align}
The adjoint system \eqref{eq:tp:adjoint}--\eqref{eq:tp:adj:TC}
is a backward-in-time parabolic problem with homogeneous terminal condition \eqref{eq:tp:adj:TC}.

\begin{rmk}
The signs of the $\alpha$-coupling in
\eqref{eq:tp:adj2}--\eqref{eq:tp:adj3} are those of the transpose of
the storage block of the state system.
\end{rmk}
Well-posedness of \eqref{eq:tp:adjoint} with \eqref{eq:tp:adj:TC} follows by time reversal and the same mechanical
Schur-complement reduction as for the state system.  After the change
$t\mapsto T-t$, the sign of the scalar coupling in the mechanical block
is reversed, so the mechanical lift of the scalar derivative changes
sign; the coupling in the scalar adjoint equation changes sign at the
same time.  Consequently the two sign changes cancel after elimination
and the reduced adjoint equation contains the same effective mass form
$m_{\rm eff}$ as \eqref{eq:tp:reduced-continuous}.  Hence
Lemma~\ref{lem:tp:reduced-mass} gives backward well-posedness.

The following lemma characterizes the optimal control via a
variational inequality and an equivalent pointwise projection formula.

\begin{lemma}
  \label{lem:tp:VI}
  Let
  $(u(\widehat{m}_p), p_T(\widehat{m}_p), p(\widehat{m}_p))$
  be the state corresponding to
  $\widehat{m}_p \in \mathcal{M}_{p,ad}$,
  and let
  $(v^*(\widehat{m}_p), q_T^*(\widehat{m}_p), q^*(\widehat{m}_p))$
  be the corresponding adjoint solution of
  \eqref{eq:tp:adjoint}--\eqref{eq:tp:adj:TC}.
  Then $\widehat{m}_p$ is the optimal control if and
  only if the following variational inequality holds:
  \begin{align}
    \label{eq:tp:VI-p}
    \int_0^T
    \bigl(
      \gamma_p \widehat{m}_p + q^*(\widehat{m}_p),
      m_p - \widehat{m}_p
    \bigr)_\Omega \, dt
    &\ge 0
    \qquad \forall m_p \in \mathcal{M}_{p,ad}.
  \end{align}
  Moreover, \eqref{eq:tp:VI-p} is equivalent to the pointwise
  projection formula
  \begin{align}
    \label{eq:tp:proj-p}
    \widehat{m}_p(t,x)
    &= \mathcal{P}_{\mathcal{M}_{p,ad}}
       \!\left(-\frac{1}{\gamma_p} q^*(\widehat{m}_p)(t,x)\right)
    \quad \text{a.e.\ in } (0,T) \times \Omega,
  \end{align}
  where the pointwise projection onto the admissible set is defined by
  \begin{align}
    \label{eq:tp:proj:def-p}
    \mathcal{P}_{\mathcal{M}_{p,ad}}(r)(t,x)
    &:= \max\{l_p, \min\{r(t,x), L_p\}\}.
  \end{align}
\end{lemma}

\begin{proof}
With the sign convention $\mathcal{L}=J-\int_0^T\mathcal{R}\,dt$,
the variation with respect to the control gives gradient
$\gamma_p m_p+q^*$.
Let $\widehat{m}_p$ be the optimal control and let
$m_p \in \mathcal{M}_{p,ad}$ be arbitrary.
Since $\mathcal{M}_{p,ad}$ is convex, for any $\delta \in (0,1]$
the perturbed control
$\widehat{m}_p + \delta(m_p - \widehat{m}_p) \in \mathcal{M}_{p,ad}$.
By optimality,
\begin{align*}
  j\bigl(\widehat{m}_p + \delta(m_p - \widehat{m}_p)\bigr)
  \ge j(\widehat{m}_p),
\end{align*}
so dividing by $\delta > 0$ and letting $\delta \to 0$ yields
$j'(\widehat{m}_p)(m_p - \widehat{m}_p) \ge 0$.

Differentiating $j$ in the direction $m_p - \widehat{m}_p$ and using
the adjoint equations to express the state sensitivities, we obtain
\begin{align*}
  j'(\widehat{m}_p)(m_p - \widehat{m}_p)
  = \int_0^T
    \bigl(
      \gamma_p \widehat{m}_p + q^*(\widehat{m}_p),
      m_p - \widehat{m}_p
    \bigr)_\Omega \,dt,
\end{align*}
which gives \eqref{eq:tp:VI-p}.
The converse direction (variational inequality implies optimality)
follows because $j$ is convex, since $J$ is convex and the
control-to-state map is affine;
see \cite{hinze2009optimization} for general background.

The variational inequality \eqref{eq:tp:VI-p} can be written
pointwise a.e.\ in $(0,T) \times \Omega$ as
\begin{align*}
  \bigl(
    \gamma_p \widehat{m}_p(t,x) + q^*(t,x),
    \eta - \widehat{m}_p(t,x)
  \bigr) \ge 0
  \qquad \forall \eta \in [l_p, L_p].
\end{align*}
This is precisely the optimality condition for the scalar projection
problem $\min_{\eta \in [l_p, L_p]}
\frac{1}{2}|\eta + q^*/\gamma_p|^2$,
whose unique solution is
$\eta = \mathcal{P}_{\mathcal{M}_{p,ad}}(-q^*/\gamma_p)$
as defined in \eqref{eq:tp:proj:def-p}.
This gives \eqref{eq:tp:proj-p}.
\end{proof}

The optimality system is coupled through the adjoint right-hand side
and the projection formula for the control.

\section{Finite element discretization and discrete optimal control}
\label{sec:THM:FEM}

\subsection{Finite element spaces}
\label{sec:THM:FE-spaces}
Let $\{\mathcal{T}_h\}_{h>0}$ be a family of shape-regular
partitions of $\overline{\Omega}$ into triangles (or tetrahedra for
$d=3$) with mesh size $h = \max_{K \in \mathcal{T}_h} h_K$.  We use
the notation $\mathbb{P}_k(K)$ for the space of polynomials of total
degree at most $k$ on $K$, and $\mathbb{P}_k^{\mathrm{dc}}(\mathcal{T}_h)$
for the space of piecewise polynomials of degree $k$ with no
interelement continuity imposed.

For the displacement and total pressure pair $(u, p_T)$ we use
the lowest-order Taylor--Hood ($\mathbb P_2/\mathbb P_1$) pair:
\begin{align}
  \label{eq:THM:VhQh0}
  V_h
  &:= \bigl\{ v_h \in V \,:\,
       v_h|_K \in [\mathbb{P}_2(K)]^d
       \;\forall\, K \in \mathcal{T}_h \bigr\},
  \\
  \label{eq:THM:Qh0}
  Q_{T,h}
  &:= \bigl\{ q_h \in H^1(\Omega) \,:\,
       q_h|_K \in \mathbb{P}_{1}(K)
       \;\forall\, K \in \mathcal{T}_h \bigr\}
  \;\subset\; H^1(\Omega) \subset L^2(\Omega) = Q_T.
\end{align}
%
It is known that the pair $(V_h, Q_{T,h})$ satisfies the discrete inf-sup condition on the full pressure space $Q_{T,h}$:
\begin{align}
  \label{eq:THM:TH-infsup}
  \inf_{0 \ne q_{T,h} \in Q_{T,h}}
  \sup_{v_h \in V_h}
  \frac{(q_{T,h}, \dive v_h)_{\Omega}}
       {\|v_h\|_1 \|q_{T,h}\|_0}
  \ge C > 0,
\end{align}
with $C$ independent of $h$; see, e.g., \cite{boffi2013mixed}.

The discrete inf-sup condition yields a right inverse
$R_h:Q_{T,h}\to V_h$ satisfying
\begin{align}
  \label{eq:THM:right-inverse}
  b(r_{T,h},R_hq_{T,h})&=-(r_{T,h},q_{T,h})_\Omega
       \qquad\forall r_{T,h},q_{T,h}\in Q_{T,h},\notag\\
  \|R_hq_{T,h}\|_V&\le C_R\|q_{T,h}\|_0,
\end{align}
where $C_R$ is independent of $h$, $\lambda$, and $s_0$.
It is the discrete
analogue of \eqref{eq:R-operator}.

For $p$ we use the standard linear Lagrange finite element:
\begin{align}
  \label{eq:THM:Qh}
  Q_{h}
  &:= \bigl\{ q_h \in H_0^1(\Omega) \,:\,
       q_h|_K \in \mathbb{P}_{1}(K)
       \;\forall\, K \in \mathcal{T}_h \bigr\} \subset Q.
\end{align}
%


\begin{lemma}
  \label{lem:tp:discrete-mech}
  For given $p\in L^2(\Omega)$, let
  $(\mathcal U_hp,\mathcal P_{T,h}p)\in V_h\times Q_{T,h}$ be the unique
  solution of
  \begin{subequations}
    \label{eq:mech-op-disc}
  \begin{align}
    \label{eq:mech-op-disc1}
    a_u(\mathcal U_hp,v_h)+b(\mathcal P_{T,h}p,v_h)&=0,
    \\
    \label{eq:mech-op-disc2}
    b(q_{T,h},\mathcal U_hp)-\lambda^{-1}(\mathcal P_{T,h}p,q_{T,h})_\Omega
      +\lambda^{-1}\alpha(p,q_{T,h})_\Omega&=0
  \end{align}
  \end{subequations}
  for all $(v_h,q_{T,h})\in V_h\times Q_{T,h}$.  Then
  $\|\mathcal U_hp\|_V+\|\mathcal P_{T,h}p\|_0\le C\|p\|_0$, with $C$
  independent of $h$ and $s_0$, and uniform for $\lambda\ge\lambda_0>0$.
  Define, for $p_h,q_h\in Q_h$,
  \begin{align}
    \label{eq:tp:discrete-mass-def}
    m_{{\rm eff},h}(p_h,q_h):=-\lambda^{-1}\alpha(\mathcal P_{T,h}p_h,q_h)_\Omega
      +s_1(p_h,q_h)_\Omega .
  \end{align}
  Then $m_{{\rm eff},h}$ is symmetric and bounded on $Q_h$, and for every $q_h\in Q_h$,
  \begin{align}
    \label{eq:tp:discrete-mass-identity}
    m_{{\rm eff},h}(q_h,q_h)
    =s_0\|q_h\|_0^2+a_u(\mathcal U_hq_h,\mathcal U_hq_h)
      +\lambda^{-1}\|\mathcal P_{T,h}q_h-\alpha q_h\|_0^2 \ge c_{{\rm eff},h}\|q_h\|_0^2
  \end{align}
  where $c_{{\rm eff},h}:=s_0 + \frac{\alpha^2}{\lambda+C_R^2}$
  with $C_R$ from \eqref{eq:THM:right-inverse}. 
\end{lemma}

\begin{proof}
  Existence, uniqueness, and the bound for
  $(\mathcal U_hp,\mathcal P_{T,h}p)$ follow from the discrete
  Babu\v{s}ka--Brezzi theory with \eqref{eq:THM:TH-infsup}, as in
  Lemma~\ref{lem:tp:mech-operator}. In the rest of proof, assume $p =q_h \in Q_h$. Then, symmetry follows by the argument of Lemma~\ref{lem:tp:reduced-mass}. Testing \eqref{eq:mech-op-disc1} with $\mathcal U_hq_h$ and \eqref{eq:mech-op-disc2} with $\mathcal P_{T,h}q_h$ gives
  \begin{align*}
    a_u(\mathcal U_h q_h,\mathcal U_h q_h)
    =\lambda^{-1}\bigl(\alpha(q_h,\mathcal P_{T,h} q_h)_{\Omega}
                       -\|\mathcal P_{T,h}q_h\|_0^2\bigr).
  \end{align*}
  The first identity in
  \eqref{eq:tp:discrete-mass-identity} follows by the argument of
  Lemma~\ref{lem:tp:reduced-mass}. 
  For the inequality in \eqref{eq:tp:discrete-mass-identity}, 
  since
  $q_h\in Q_h\subset Q_{T,h}$, $R_h q_h$ is well-defined. Testing the first equation of
  \eqref{eq:mech-op-disc} with $v_h=R_hq_h$ and using
  \eqref{eq:THM:right-inverse} with $r_{T,h}=\mathcal{P}_{T,h}q_h$ gives
  $a_u(\mathcal{U}_h q_h,R_hq_h)-(\mathcal{P}_{T,h}q_h,q_h)_\Omega=0$. From this
  \begin{align*}
    \alpha\|q_h\|_0^2
    &= (\alpha q_h - \mathcal{P}_{T,h} q_h, q_h)_{\Omega} + a_u(\mathcal{U}_h,R_hq_h) 
    \\
    &\le (\|\alpha q_h - \mathcal{P}_{T,h}q_h\|_0 + C_R\|\mathcal{U}_h q_h\|_V ) \|q_h\|_0
    \\
    &\le (\lambda +C_R^2)^{\frac 12} \left( \lambda^{-1}\|\alpha q_h - \mathcal{P}_{T,h}q_h\|_0^2 + \|\mathcal{U}_h q_h\|_V^2\right)^{\frac 12} \|q_h \|_0 .
  \end{align*}
  Combining it with the identity in \eqref{eq:tp:discrete-mass-identity}, we get the inequality in \eqref{eq:tp:discrete-mass-identity}.
\end{proof}

\begin{rmk}
  \label{rmk:THM:lambda-uniform-FD}
  The proof of Lemma~\ref{lem:tp:discrete-mech} relies on $Q_h \subset Q_{T,h}$. 
\end{rmk}

\subsection{Time discretization and fully discrete state system}
\label{sec:THM:time-disc}

Let $\Delta t = T/N$ be the time step size with positive integer $N$,
and set $t_k = k\Delta t$ for $k = 0, 1, \ldots, N$.
We write $I_k:=(t_{k-1},t_k]$, $k=1,\ldots,N$.  For a time-dependent
function $z$ we denote its average over $I_k$ by
\begin{align}
  \label{eq:tp:cell-avg}
  \bar z^k := \frac{1}{\Delta t}\int_{I_k} z(t)\,dt .
\end{align}
Superscripts without an overbar are reserved for discrete time levels
(e.g.\ $u_h^k$).
For the control variable we use the space-time piecewise
constant space
\begin{align}
  \label{eq:THM:Mh}
  M_{p,h}
    &:= \bigl\{ m \in L^2(0,T;L^2(\Omega)) \,:\,
       m|_{I_k\times K} \in \mathbb{P}_0(I_k)\otimes\mathbb{P}_0(K)
       \;\forall\, k=1,\ldots,N,\; K \in \mathcal{T}_h \bigr\},
\end{align}
i.e.\ functions that are piecewise constant on each space-time cell
$I_k\times K$.
For a sequence $\{g^k\}_{k \ge 0}$ we use the backward difference
quotient
\begin{align*}
  d_t g^{k+1} := \frac{g^{k+1} - g^k}{\Delta t}.
\end{align*}
We treat $\{m_{p,h}^k\}_{k=1}^N$ as
space-time functions that are constant on each interval $(t_{k-1},t_k]$
with values $m_{p,h}^k$, and similarly for
the state and adjoint sequences.
For a sequence $\{g^k\}_{k=1}^N$ in a normed space $X$ we write
\begin{align*}
  \|\{g^k\}\|_{\ell^2_tX}^2
  := \Delta t \sum_{k=1}^N \|g^k\|_X^2 ;
\end{align*}
for piecewise constant-in-time functions this coincides with the
$L^2(0,T;X)$ norm.

Suppose that the numerical initial data satisfy
\begin{align*}
  \|u(0) - u_h^0\|_1
  + \|p_T(0) - p_{T,h}^0\|_0
  + \|p(0) - p_h^0\|_0
  \le C_{\rm init}\,h,
\end{align*}
with $C_{\rm init}$ independent of $h$ and $\Delta t$, and of
$\lambda\ge\lambda_0$ for the $\lambda$-uniform estimates.
Given $m_p\in\mathcal M_{p,ad}$, the fully discrete
state solution at the $(k+1)$-th time step is defined by backward Euler:
find
$(u_h^{k+1}, p_{T,h}^{k+1}, p_h^{k+1})
\in V_h \times Q_{T,h} \times Q_h$
such that
\begin{subequations}
  \label{eq:tp:FD-state}
  \begin{align}
    \label{eq:tp:FD-state1}
    a_u(d_t u_h^{k+1}, v_h)
    + b(d_t p_{T,h}^{k+1}, v_h)
    &= \frac{1}{\Delta t}\int_{t_k}^{t_{k+1}} (f, v_h)_\Omega \,dt,
    \\
    \label{eq:tp:FD-state2}
    b(q_{T,h}, d_t u_h^{k+1})
    - \lambda^{-1}(d_t p_{T,h}^{k+1}, q_{T,h})_\Omega
    + \lambda^{-1}\alpha (d_t p_h^{k+1}, q_{T,h})_\Omega
    &= 0,
    \\
    \label{eq:tp:FD-state3}
    -\lambda^{-1}\alpha (d_t p_{T,h}^{k+1}, q_h)_\Omega
    + s_1(d_t p_h^{k+1}, q_h)_\Omega
    + (\kappa_p \nabla p_h^{k+1}, \nabla q_h)_\Omega
    &= \frac{1}{\Delta t}\int_{t_k}^{t_{k+1}} (m_p, q_h)_\Omega \,dt,
  \end{align}
\end{subequations}
for all $(v_h, q_{T,h}, q_h) \in V_h \times
Q_{T,h} \times Q_h$.
For the discrete optimal control problem, $m_p=m_{p,h}\in\mathcal M_{p,ad}^h$,
and the right-hand side of \eqref{eq:tp:FD-state3} is
$(m_{p,h}^{k+1},q_h)_\Omega$.

\begin{lemma}[discrete forward stability]
  \label{lem:tp:generic-stability}
  Let $(w_h^l,\xi_h^l,r_h^l)_{l=0}^N\subset V_h\times Q_{T,h}\times Q_h$
  satisfy, for $l=1,\ldots,N$,
  \begin{subequations}
    \label{eq:tp:generic-system}
    \begin{align}
      a_u(d_tw_h^l,v_h)+b(d_t\xi_h^l,v_h)
      &=(F_w^l,v_h)_\Omega,\label{eq:tp:gen1}\\
      b(q_{T,h},d_tw_h^l)-\lambda^{-1}(d_t\xi_h^l,q_{T,h})_\Omega
      +\lambda^{-1}\alpha(d_tr_h^l,q_{T,h})_\Omega
      &=0,\label{eq:tp:gen2}\\
      -\lambda^{-1}\alpha(d_t\xi_h^l,q_h)_\Omega
      +s_1(d_tr_h^l,q_h)_\Omega+(\kappa_p\nabla r_h^l,\nabla q_h)_\Omega
      &=(F_r^l,q_h)_\Omega.\label{eq:tp:gen3}
    \end{align}
  \end{subequations}
  Assume $s_0\ge0$ and zero initial data.
  Then
  \begin{align}
    \label{eq:tp:generic-stability}
    \max_{1\le l\le N}\bigl(
      \|w_h^l\|_V+\|\xi_h^l\|_0+\|r_h^l\|_0\bigr)
    \le C_{\rm eff}\Bigl(\Delta t\sum_{l=1}^N
      (\|F_w^l\|_{V^*}^2+\|F_r^l\|_0^2)\Bigr)^{1/2},
  \end{align}
  where $C_{\rm eff}$ is independent of $h$, $\Delta t$, but depends on  $c_{{\rm eff},h}^{-1}$
  in \eqref{eq:tp:discrete-mass-identity}.
\end{lemma}

\begin{proof}
  At each time level, split the differentiated mechanical variables as
  \[
    d_tw_h^l=w_{F,h}^l+\mathcal U_h(d_tr_h^l),
    \qquad
    d_t\xi_h^l=\xi_{F,h}^l+\mathcal P_{T,h}(d_tr_h^l),
  \]
  where $(w_{F,h}^l,\xi_{F,h}^l)$ solves the discrete mechanical
  problem with right-hand side $F_w^l$ and zero scalar datum.  The
  discrete saddle-point estimate gives
  \[
    \|w_{F,h}^l\|_V+\|\xi_{F,h}^l\|_0
    \le C\|F_w^l\|_{V^*},
  \]
  uniformly in $h$.  Substituting into \eqref{eq:tp:gen3} gives
  \[
    m_{{\rm eff},h}(d_tr_h^l,q_h)+(\kappa_p\nabla r_h^l,\nabla q_h)_\Omega
    =(F_r^l,q_h)_\Omega
      +\lambda^{-1}\alpha(\xi_{F,h}^l,q_h)_\Omega .
  \]
  Taking $q_h=d_tr_h^l$, using Lemma~\ref{lem:tp:discrete-mech},
  multiplying by $\Delta t$, summing, and applying Young's inequality
  yield an $\ell_t^2L_x^2$ bound for $d_tr_h$.  The mechanical decomposition and
  Lemma~\ref{lem:tp:discrete-mech} give the corresponding bounds for
  $d_tw_h$ and $d_t\xi_h$.  Since the initial data vanish, discrete
  summation in time and Cauchy--Schwarz give
  \eqref{eq:tp:generic-stability}.
\end{proof}

\begin{rmk}
  The state system and the state-difference system have exactly the
  sign pattern of \eqref{eq:tp:generic-system}.  The
  adjoint system, after time reversal, has the same stability
  structure as the forward state scheme.
\end{rmk}

\begin{lemma}[Well-posedness and stability of the discrete state system]
  \label{lem:tp:FD-wellposed}
  The fully discrete state system \eqref{eq:tp:FD-state} is
  well-posed at each time step, and its solution with zero initial
  data satisfies
  \begin{align}
    \label{eq:tp:FD-state-stability}
    \max_{1\le k\le N}\bigl(
      \|u_h^k\|_V + \|p_{T,h}^k\|_0
      + \|p_h^k\|_0\bigr)
    \le C\Bigl(\Delta t\sum_{k=1}^N\bigl(
      \|f_u^k\|_{V^*}^2+\|f_p^k\|_0^2
    \bigr)\Bigr)^{1/2},
  \end{align}
  where $f_u^k:=\frac{1}{\Delta t}\int_{t_{k-1}}^{t_k}f\,dt$ and
  $f_p^k:=\frac{1}{\Delta t}\int_{t_{k-1}}^{t_k}m_p\,dt$ are the
  right-hand sides of \eqref{eq:tp:FD-state1} and \eqref{eq:tp:FD-state3}
  at step $k$,
  with $C$ independent of $h$, $\Delta t$.
\end{lemma}

\begin{proof}
    It is straightforward from Lemma~\ref{lem:tp:generic-stability}, so we omit details.
\end{proof}

\subsection{Discrete optimal control problem}
\label{sec:THM:discrete-OCP}

We define the discrete admissible set
\begin{align*}
  \mathcal{M}_{p,ad}^h
  &:= \mathcal{M}_{p,ad} \cap M_{p,h}
   = \bigl\{ m_h \in M_{p,h} : l_p \le m_h \le L_p \ \text{a.e.} \bigr\}.
\end{align*}
Note that $\mathcal{M}_{p,ad}^h \subset \mathcal{M}_{p,ad}$ since
$M_{p,h}$ is a piecewise constant subspace of $L_t^2 L_x^2$.

We also introduce the $L^2((0,T)\times\Omega)$-orthogonal projection
$\Pi_h^M : L_t^2L_x^2 \to M_{p,h}$,
defined cellwise in space and time by
\begin{align}
  \label{eq:tp:PiM-def}
  (\Pi_h^M m)\big|_{(t_{k-1},t_k)\times K}
  := \frac{1}{\Delta t\,|K|}\int_{t_{k-1}}^{t_k}\!\!\int_K m\,dx\,dt,
  \qquad K \in \mathcal{T}_h,\ k=1,\ldots,N.
\end{align}
Since $\Pi_h^M m$ is a cellwise space-time average of $m$, if
$l_p \le m \le L_p$ a.e.\ then $l_p \le \Pi_h^M m \le L_p$ a.e.;
hence $\Pi_h^M$ maps $\mathcal{M}_{p,ad}$ into $\mathcal{M}_{p,ad}^h$:
\begin{align}
  \label{eq:tp:PiM-admissible}
  \Pi_h^M : \mathcal{M}_{p,ad} \to \mathcal{M}_{p,ad}^h.
\end{align}

The discrete objective functional is
\begin{align}
  \label{eq:tp:Jh}
  &J_h\bigl(\{u_h^k\}_{k=1}^N, \{p_h^k\}_{k=1}^N,
    \{m_{p,h}^k\}_{k=1}^N\bigr)
  \notag \\
  &:=
  \frac{\alpha_{u,C}}{2} \Delta t \sum_{k=0}^{N-1}
    \|u_h^{k+1} - \bar u_C^{k+1}\|_0^2
  + \frac{\alpha_{p,C}}{2} \Delta t \sum_{k=0}^{N-1}
    \|p_h^{k+1} - \bar p_C^{k+1}\|_0^2
  + \frac{\gamma_p}{2} \|\{m_{p,h}^k\}\|_{\ell^2_tL^2_x}^2,
\end{align}
subject to the fully discrete state system \eqref{eq:tp:FD-state},
where $(u_h^N, p_h^N)$ is the state at the final time
step $t_N = T$.  As in the continuous setting, we use the discrete
reduced functional $j_h(m_{p,h}):=J_h(\{u_h^k\},\{p_h^k\},\{m_{p,h}^k\})$
on $\mathcal M_{p,ad}^h$, where $(u_h^k,p_{T,h}^k,p_h^k)$ solves
\eqref{eq:tp:FD-state} with control $m_{p,h}$.

\begin{theorem}
  \label{thm:tp:discrete-existence}
  There exists a unique solution
  $\widehat{m}_{p,h} \in \mathcal{M}_{p,ad}^h$
  to the discrete optimal control problem
  \eqref{eq:tp:FD-state}--\eqref{eq:tp:Jh}.
\end{theorem}

\begin{proof}
  The proof follows the same lines as the continuous case in
  Theorem~\ref{thm:tp:OCP-existence}, so we omit here.
\end{proof}

\subsection{Discrete first-order optimality system}
\label{sec:THM:discrete-optimality}

\begin{theorem}
  \label{thm:tp:discrete-optimality}
  A control $\widehat{m}_{p,h} \in \mathcal{M}_{p,ad}^h$
  and the corresponding state
  $\{(\widehat{u}_h^{k+1}, \widehat{p}_{T,h}^{k+1}, \widehat{p}_h^{k+1})\}_{k=0}^{N-1}$
  is the optimal solution to
  \eqref{eq:tp:FD-state}--\eqref{eq:tp:Jh} if and only if there
  exists an adjoint state
  $\{(\widehat{v}_h^{*,k}, \widehat{q}_{T,h}^{*,k}, \widehat{q}_h^{*,k})\}_{k=0}^{N-1}$
  satisfying the discrete adjoint equations
  \begin{subequations}
    \label{eq:tp:FD-adj}
    \begin{align}
      \label{eq:tp:FD-adj1}
      -a_u(d_t v_h^{*,k+1}, v_h)
      - b(d_t q_{T,h}^{*,k+1}, v_h)
      &= \alpha_{u,C}
         (u_h^{k+1} - \bar u_C^{k+1}, v_h)_\Omega,
      \\
      \label{eq:tp:FD-adj2}
      -b(q_{T,h}, d_t v_h^{*,k+1})
      + \lambda^{-1}(d_t q_{T,h}^{*,k+1}, q_{T,h})_\Omega
      + \lambda^{-1}\alpha (d_t q_h^{*,k+1}, q_{T,h})_\Omega
      &= 0,
      \\
      \label{eq:tp:FD-adj3}
      -\lambda^{-1}\alpha (d_t q_{T,h}^{*,k+1}, q_h)_\Omega
      - s_1(d_t q_h^{*,k+1}, q_h)_\Omega
      + (\kappa_p \nabla q_h^{*,k}, \nabla q_h)_\Omega
      &= \alpha_{p,C}
         (p_h^{k+1} - \bar p_C^{k+1}, q_h)_\Omega,
    \end{align}
  \end{subequations}
  for all $(v_h, q_{T,h}, q_h) \in V_h
  \times Q_{T,h} \times Q_h$ and $k = 0, 1, \ldots, N-1$,
  where $\bar u_C^{k+1}$, $\bar p_C^{k+1}$ denote the
  cell-average values of the target fields on the interval
  $(t_k,t_{k+1})$.
  \begin{align}
    \label{eq:tp:FD-adj-TC}
    v_h^{*,N} = 0,\quad
    q_{T,h}^{*,N} = 0,\quad
    q_h^{*,N} = 0.
  \end{align}

    Together with the discrete variational inequality
  \begin{align}
    \label{eq:tp:FD-VI-p}
    \Delta t \sum_{k=0}^{N-1}
    \bigl(\gamma_p \widehat{m}_{p,h}^{k+1} + q_h^{*,k}(\widehat{m}_{p,h}),\,
    z_p^{k+1} - \widehat{m}_{p,h}^{k+1}\bigr)_\Omega
    &\ge 0
    \qquad \forall\, z_p \in \mathcal{M}_{p,ad}^h.
  \end{align}
  Moreover, \eqref{eq:tp:FD-VI-p} is equivalent to the cellwise
  projection formula: for $k = 1, 2, \ldots, N$ and $K\in\mathcal T_h$,
  \begin{align}
    \label{eq:tp:FD-proj-p}
    \widehat{m}_{p,h}^k\big|_K
    &= \max\Bigl\{l_p,\ \min\Bigl\{-\frac{1}{\gamma_p|K|}
       \int_K q_h^{*,k-1}(\widehat{m}_{p,h})\,dx,\ L_p\Bigr\}\Bigr\}.
  \end{align}
\end{theorem}

\begin{proof}
  The proof follows the same lines as the continuous
  Lemma~\ref{lem:tp:VI}.
  By the optimality of $\widehat{m}_{p,h}$ and the
  convexity of $\mathcal{M}_{p,ad}^h$, for any admissible
  perturbation direction and $\delta \in (0,1]$ the G\^ateaux
  derivative of $j_h$ satisfies $j_h'(\widehat{m}_{p,h})(z_p - \widehat{m}_{p,h}) \ge 0$.
  Differentiating $j_h$ using the discrete state equations
  \eqref{eq:tp:FD-state} and the discrete adjoint equations
  \eqref{eq:tp:FD-adj}--\eqref{eq:tp:FD-adj-TC} to represent state
  sensitivities gives \eqref{eq:tp:FD-VI-p}.
  Testing \eqref{eq:tp:FD-VI-p} with $z_p$ equal to $\widehat m_{p,h}$
  except on one space-time cell $I_k\times K$, where $z_p=\eta\in[l_p,L_p]$,
  and dividing by $\Delta t|K|>0$ gives the scalar variational inequality
  $\bigl(\gamma_p\widehat m_{p,h}^k|_K+|K|^{-1}\int_K q_h^{*,k-1}\,dx\bigr)
  \bigl(\eta-\widehat m_{p,h}^k|_K\bigr)\ge0$ for all
  $\eta\in[l_p,L_p]$, whose unique solution is \eqref{eq:tp:FD-proj-p}.
  Conversely, \eqref{eq:tp:FD-proj-p} gives these cellwise inequalities,
  and summing them over all cells gives \eqref{eq:tp:FD-VI-p}.
\end{proof}

\begin{rmk}
  \label{rmk:THM:discrete-adjoint-indexing}
  The control value on $(t_{k-1},t_k]$ is paired with
  $q_h^{*,k-1}$.  The terminal value
  $q_h^{*,N}=0$ serves only as starting datum for the
  backward adjoint recursion.
\end{rmk}

We record a stability estimate for the discrete adjoint
system with general interior right-hand sides and homogeneous terminal
data, which is the discrete analogue of the backward parabolic energy
estimate.  It is used in the sensitivity and adjoint error estimates of
Section~\ref{sec:THM:error}.

\begin{lemma}[Discrete adjoint stability]
  \label{lem:tp:FD-adj-stability}
  Let $(v_h^{*,k},q_{T,h}^{*,k},q_h^{*,k})_{k=0}^{N-1}$ solve the
  discrete adjoint recursion \eqref{eq:tp:FD-adj} with general interior
  data $\alpha_{u,C}F_u^{k+1}$ and $\alpha_{p,C}F_p^{k+1}$ and with
  homogeneous terminal condition.  Then
  \begin{align}
    \label{eq:tp:FD-adj-stability}
    \max_{0\le k\le N-1}\bigl(
      \|v_h^{*,k}\|_V+\|q_{T,h}^{*,k}\|_0+\|q_h^{*,k}\|_0\bigr)
    \le C_{\rm eff}\Bigl(\Delta t\sum_{j=1}^{N}
      (\|F_u^j\|_0^2+\|F_p^j\|_0^2)\Bigr)^{1/2},
  \end{align}
  where $C_{\rm eff}$ is independent of $h$, $\Delta t$, but may depend on the lower bound of $c_{{\rm eff},h}$ 
\end{lemma}

\begin{proof}
  Reverse time by setting
  $(\check v_h^l,\check q_{T,h}^l,\check q_h^l)
   =(v_h^{*,N-l},q_{T,h}^{*,N-l},q_h^{*,N-l})$.
  The reversed equations are, for $l=1,\ldots,N$,
  \begin{subequations}
    \label{eq:tp:adj-reversed}
    \begin{align}
      a_u(d_t\check v_h^l,v_h)+b(d_t\check q_{T,h}^l,v_h)
      &=\alpha_{u,C}(F_u^{N-l+1},v_h)_\Omega,\\
      b(q_{T,h},d_t\check v_h^l)-\lambda^{-1}(d_t\check q_{T,h}^l,q_{T,h})_\Omega
      -\lambda^{-1}\alpha(d_t\check q_h^l,q_{T,h})_\Omega&=0,\\
      \lambda^{-1}\alpha(d_t\check q_{T,h}^l,q_h)_\Omega
      +s_1(d_t\check q_h^l,q_h)_\Omega+(\kappa_p\nabla\check q_h^l,\nabla q_h)_\Omega
      &=\alpha_{p,C}(F_p^{N-l+1},q_h)_\Omega.
    \end{align}
  \end{subequations}
  Their mechanical block differs from the forward state only in the sign of
  the scalar coupling.  Consequently,
  \[
    d_t\check v_h^l=w_{F,h}^l-\mathcal U_h(d_t\check q_h^l),
    \qquad
    d_t\check q_{T,h}^l=\xi_{F,h}^l-\mathcal P_{T,h}(d_t\check q_h^l),
  \]
  where $(w_{F,h}^l,\xi_{F,h}^l)$ is the mechanical lift of the
  reversed displacement datum.  Substitution in the third reversed
  adjoint equation gives
  \[
    m_{{\rm eff},h}(d_t\check q_h^l,q_h)+(\kappa_p\nabla\check q_h^l,\nabla q_h)_\Omega
    =\alpha_{p,C}(F_p^{N-l+1},q_h)_\Omega
      -\lambda^{-1}\alpha(\xi_{F,h}^l,q_h)_\Omega .
  \]
  Thus the same effective mass $m_{{\rm eff},h}$ occurs as in the forward
  reduction.  Testing with $q_h=d_t\check q_h^l$, summing in time, and
  using Lemma~\ref{lem:tp:discrete-mech} gives the scalar derivative
  bound without $s_0^{-1}$.  Mechanical reconstruction gives the
  remaining derivative bounds.  Since the reversed initial data are
  zero, discrete integration in time yields
  \eqref{eq:tp:FD-adj-stability}; reversing time completes the proof.
\end{proof}

\section{A priori error analysis}
\label{sec:THM:error}

The well-posedness and stability results above are valid for $s_0\ge0$
with constants independent of $s_0^{-1}$, at the price that the
effective-mass coercivity constant may depend on $\lambda$.  In this
section, when a $\lambda$-uniform error constant is asserted, we use the
alternative bound $m_{\rm eff}(p,p)\ge s_0\|p\|_0^2$ and assume that $s_0$ is
fixed positive.  Thus storage-robustness and
near-incompressibility robustness are distinguished explicitly.

In this section we derive a priori error estimates for the fully
discrete scheme \eqref{eq:tp:FD-state}--\eqref{eq:tp:FD-proj-p}.
Throughout this section, we assume that the family of triangulations is shape-regular and quasi-uniform.

\textbf{Regularity assumption (Assumption~R).}
Throughout Section~\ref{sec:THM:error}
we assume that the continuous optimal state and adjoint satisfy
\begin{align}
  \notag
  \widehat u,\,\widehat v^*
  &\in L^2(0,T;\boldsymbol{H}^2(\Omega))\cap H^1(0,T;V),
  \\
  \label{eq:tp:reg-state-adj}
  \widehat p_T,\,\widehat q_T^*
  &\in L^2(0,T;H^1(\Omega))\cap H^1(0,T;Q_T),
  \\
  \notag
  \widehat p,\,\widehat q^{\,*}
  &\in L^2(0,T;H^2(\Omega))\cap H^1(0,T;L^2(\Omega)),
\end{align}
with all constants uniform for $\lambda\ge\lambda_0>0$.
Sufficient conditions for \eqref{eq:tp:reg-state-adj} are established
in Appendix~\ref{app:regularity}. The control regularity used below
follows from the projection formula, as shown in
Remark~\ref{rmk:tp:control-regularity}.

For the continuous optimal solution we write
$(\widehat{u}, \widehat{p}_T, \widehat{p},
 \widehat{m}_p,
 \widehat{v}^*, \widehat{q}_T^*, \widehat{q}^*)$,
and for the discrete optimal solution
$\{(\widehat{u}_h^{k+1}, \widehat{p}_{T,h}^{k+1}, \widehat{p}_h^{k+1},
   \widehat{m}_{p,h}^{k+1},
   \widehat{v}_h^{*,k}, \widehat{q}_{T,h}^{*,k}, \widehat{q}_h^{*,k})\}$.

\subsection{Auxiliary discrete systems}
\label{sec:THM:auxiliary}

For any $m_p \in \mathcal{M}_{p,ad}$ we define
the auxiliary fully discrete state solution
$(u_h^{k+1}(m_p), p_{T,h}^{k+1}(m_p),
 p_h^{k+1}(m_p))$
by \eqref{eq:tp:FD-state} with right-hand side
$\frac{1}{\Delta t}\int_{t_k}^{t_{k+1}} (m_p, q_h)_\Omega \,dt$.
We also define the auxiliary discrete adjoint solution
$(v_h^{*,k}(m_p), q_{T,h}^{*,k}(m_p),
 q_h^{*,k}(m_p))$
by \eqref{eq:tp:FD-adj}--\eqref{eq:tp:FD-adj-TC}
driven by the auxiliary state.

For $m_{p,h}\in\mathcal M_{p,ad}^h$, let $q_h^*(m_{p,h})$ denote the
piecewise constant-in-time function equal to $q_h^{*,k}(m_{p,h})$ on
$I_{k+1}$.  The $L^2_tL^2_x$-gradient of $j_h$ in $M_{p,h}$ is
\begin{align}
  \label{eq:tp:Jh-gradient}
  j_h'(m_{p,h})
  = \gamma_p m_{p,h} + \Pi_h^M q_h^*(m_{p,h}),
\end{align}
that is, $j_h'(m_{p,h})[\eta_h]
=(\gamma_pm_{p,h}+q_h^*(m_{p,h}),\eta_h)_{L^2_tL^2_x}$ for all
$\eta_h\in M_{p,h}$.

\subsection{Sensitivity of the discrete state and adjoint to the control}
\label{sec:THM:contraction}

Let $\Phi_{\rm fwd}$ denote the map
$m_{p} \mapsto (u_h, p_{T,h}, p_h)$ defined by
\eqref{eq:tp:FD-state}, and $\Phi_{\rm adj}$ the adjoint map
\eqref{eq:tp:FD-adj}.  The following lemma quantifies the sensitivity
of $\Phi_{\rm fwd}$ and $\Phi_{\rm adj}\circ\Phi_{\rm fwd}$ to the
control.

\begin{lemma}[Control-to-state and control-to-adjoint sensitivity]
  \label{lem:tp:control-sensitivity}
  Let $m_p,\tilde m_p \in
  \mathcal{M}_{p,ad}$ be two admissible controls (continuous or
  discrete; $\mathcal{M}_{p,ad}^h\subset\mathcal{M}_{p,ad}$ is
  the special case of piecewise-constant controls), and let the
  corresponding auxiliary fully discrete state and adjoint solutions
  be $\Phi_{\rm fwd}(m_p)$, $\Phi_{\rm fwd}(\tilde m_p)$
  and $\Phi_{\rm adj}\circ\Phi_{\rm fwd}$ applied to
  each.  Define the control difference $m_{p,d}:=m_p-\tilde m_p$ and
  its cell average
  \begin{align*}
    \bar m_{p,d}^{k+1} := \frac{1}{\Delta t}\int_{t_k}^{t_{k+1}}
    m_{p,d}(t)\,dt,
  \end{align*}
  which is equal to
  the value of $m_{p,d}$ on $I_{k+1}$
  when the controls are already piecewise constant in time.  By
  Jensen's inequality,
  \begin{align}
    \label{eq:tp:Md-Jensen}
    \Delta t\sum_{k=0}^{N-1}\|\bar m_{p,d}^{k+1}\|_0^2
    \le \|m_{p,d}\|_{L^2_tL^2_x}^2.
  \end{align}
  Let $(u_{d,h}^k, p_{T,d,h}^k,
  p_{d,h}^k) :=
  \Phi_{\rm fwd}(m_p) - \Phi_{\rm fwd}(\tilde m_p)$ and let
  $(v_{d,h}^{*,k}, q_{T,d,h}^{*,k}, q_{d,h}^{*,k})$
  be the corresponding difference of the adjoint solutions
  $\Phi_{\rm adj}\circ\Phi_{\rm fwd}$, with homogeneous terminal data.
  Then
  \begin{align}
    \label{eq:tp:diff-state-max}
    \max_{1\le k\le N}\bigl(
      \|u_{d,h}^k\|_V + \|p_{T,d,h}^k\|_0
      + \|p_{d,h}^k\|_0\bigr)
    &\lesssim \Bigl(\Delta t\sum_{k=1}^N
    \|\bar m_{p,d}^k\|_0^2\Bigr)^{1/2},
    \\
    \label{eq:tp:diff-state-est}
    \Delta t \sum_{k=1}^N \bigl[
      a_u(u_{d,h}^k, u_{d,h}^k) + \|p_{T,d,h}^k\|_0^2
      + \|p_{d,h}^k\|_0^2
    \bigr]
    &\lesssim C_T \Delta t \sum_{k=1}^N
    \|\bar m_{p,d}^k\|_0^2,
    \\
    \label{eq:tp:diff-adj-est}
    \Delta t \sum_{k=0}^{N-1}
    \|q_{d,h}^{*,k}\|_0^2
    &\lesssim C_T \Delta t \sum_{k=1}^N
    \|\bar m_{p,d}^k\|_0^2,
  \end{align}
  with constants independent of $h$, $\Delta t$, but may depend on $c_{{\rm eff},h}^{-1}$.
\end{lemma}

\begin{proof}
  Subtracting the auxiliary state systems for $m_p$ and $\tilde m_p$,
  which have the same initial data and body force, the state difference
  satisfies
  \begin{subequations}
    \label{eq:tp:diff-fwd}
    \begin{align}
      a_u(d_t u_{d,h}^{k+1}, v_h)
      + b(d_t p_{T,d,h}^{k+1}, v_h) &= 0,
      \\
      b(q_{T,h}, d_t u_{d,h}^{k+1})
      - \lambda^{-1}(d_t p_{T,d,h}^{k+1}, q_{T,h})_\Omega
      + \lambda^{-1}\alpha(d_t p_{d,h}^{k+1}, q_{T,h})_\Omega &= 0,
      \\
      -\lambda^{-1}\alpha(d_t p_{T,d,h}^{k+1}, q_h)_\Omega
      + s_1(d_t p_{d,h}^{k+1}, q_h)_\Omega
      + (\kappa_p \nabla p_{d,h}^{k+1}, \nabla q_h)_\Omega
      &= (\bar m_{p,d}^{k+1}, q_h)_\Omega,
    \end{align}
  \end{subequations}
  with zero initial data.  This is \eqref{eq:tp:generic-system} with
  $F_w^{k+1}=0$ and $F_r^{k+1}=\bar m_{p,d}^{k+1}$, so
  Lemma~\ref{lem:tp:generic-stability} gives \eqref{eq:tp:diff-state-max}.
  Since $\Delta t\sum_{k=1}^N a_k\le T\max_{1\le k\le N}a_k$ for
  nonnegative $a_k$, squaring \eqref{eq:tp:diff-state-max} gives
  \eqref{eq:tp:diff-state-est}.  The adjoint difference solves
  \eqref{eq:tp:FD-adj} with interior data $F_u^{k+1}=u_{d,h}^{k+1}$,
  $F_p^{k+1}=p_{d,h}^{k+1}$ and homogeneous terminal data.
  Lemma~\ref{lem:tp:FD-adj-stability}, together with
  $\|u_{d,h}^{k}\|_0\le C\|u_{d,h}^{k}\|_V$ and
  \eqref{eq:tp:diff-state-est}, bounds
  $\max_{0\le k\le N-1}\|q_{d,h}^{*,k}\|_0$ by
  $C(\Delta t\sum_k\|\bar m_{p,d}^k\|_0^2)^{1/2}$; the same summation
  inequality yields \eqref{eq:tp:diff-adj-est}.
\end{proof}

\subsection{Notation and projection operators}
\label{sec:THM:splitting}

We use the cell averages \eqref{eq:tp:cell-avg}.  The hat denotes optimal quantities: $\widehat u$ is the continuous optimal
state associated with the optimal control $\widehat m_p$, and
$\widehat u_h^j$ is the discrete optimal state at time level $j$.  Thus
$\overline{\widehat u}{}^{\,j}$, $\overline{\widehat p}{}^{\,j}$, and
$\overline{\widehat q^{\,*}}{}^{\,k+1}$ are cell averages of continuous
optimal quantities, not discrete quantities.  The error estimate in
Theorem~\ref{thm:tp:state-error} is stated in terms of these cell averages.

We use the following standard interpolation and projection operators.

\begin{definition}
  \label{def:tp:interp} Define the following interpolation operators.
  \begin{enumerate}
  \item $\Pi_h^u : \boldsymbol{H}^2(\Omega)\cap V \to
    V_h$ is the standard Scott--Zhang (or
    Cl\'ement-type) quasi-interpolation operator onto the
    Taylor--Hood velocity space, satisfying
    $\|v - \Pi_h^u v\|_1 \le Ch|v|_2$ for $v\in\boldsymbol{H}^2(\Omega)\cap V$.
  \item $\Pi_h^{p_T} : L^2(\Omega) \to Q_{T,h}$ is the
    $L^2(\Omega)$-orthogonal projection onto the pressure space
    $Q_{T,h}$, satisfying $\|p_T - \Pi_h^{p_T}p_T\|_0 \le
    Ch|p_T|_1$ for $p_T\in H^1(\Omega)$.
  \item $P_h : L^2(\Omega) \to Q_h$ is the
    $L^2(\Omega)$-orthogonal projection onto the
    diffusion space, defined by
    \[
      (P_h z, q_h)_{L^2} = (z, q_h)_{L^2}
      \quad\forall q_h\in Q_h.
    \]
  \end{enumerate}
\end{definition}
  \begin{lemma}      
    \label{lem:interpolation-property}
    The following properties hold:
    \begin{enumerate}[(a)]
      \item \emph{Orthogonality:} $((I-P_h)z,q_h)_{L^2}=0$ for all
        $q_h\in Q_h$.
      \item \emph{$L^2$-approximation:}
        $\|(I-P_h)z\|_0\le Ch\|z\|_{H^1}$ for $z\in H^1_0(\Omega)$.
      \item \emph{$H^1$-approximation:} 
        $\|(I-P_h)z\|_{H^1}\le Ch\|z\|_{H^2}$ for
        $z\in H^2(\Omega)\cap H^1_0(\Omega)$
        (see, e.g.,
        Crouzeix--Thom\'ee~\cite{CrouzeixThomee1987}).
    \end{enumerate}
  \end{lemma}

\begin{lemma}[Fixed-control state discretization error]
  \label{lem:tp:fixed-control-state-error}
  Let $(u,p_T,p)$ be the continuous state associated with
  control $m_p\in\mathcal{M}_{p,ad}$, and let
  $(u_h^j,p_{T,h}^j,p_h^j)$ be the fully discrete state
  \eqref{eq:tp:FD-state} driven by the cell averages
  $\bar m_p^{k+1}$.  Assume the corresponding state satisfies the regularity conditions in \eqref{eq:tp:reg-state-adj}.  Then
  \begin{align}
    \label{eq:tp:fixed-state-diff}
    \Bigl(\Delta t\sum_{j=1}^N
    \|\bar p^j - p_h^j\|_0^2
    \Bigr)^{1/2}
    &\le C(h + \Delta t),
    \\
    \label{eq:tp:fixed-state-mech}
    \Bigl(\Delta t\sum_{j=1}^N
    \bigl(\|\bar u^j - u_h^j\|_1^2 + \|\bar p_T^j - p_{T,h}^j\|_0^2\bigr)
    \Bigr)^{1/2}
    &\le C(h + \Delta t),
  \end{align}
  where $C>0$ is independent of $h$ and $\Delta t$ but may depends on $c_{\rm eff}^{-1}$, $c_{{\rm eff},h}^{-1}$, and $\kappa_p$.
\end{lemma}

The proof is organized through a sequence of auxiliary lemmas
(Lemmas~\ref{lem:tp:scheme-equiv}--\ref{lem:tp:mech-recovery}).
For the proof we use the continuous and discrete mechanical Schur complements introduced in Lemmas \ref{lem:tp:mech-operator} and \ref{lem:tp:discrete-mech}.

For each time slab $I_{k+1}=(t_k,t_{k+1}]$, set
\[
  \bar f^{k+1}:=\frac{1}{\Delta t}\int_{I_{k+1}}f(t)\,dt.
\]
Let $(u_{f,h}^{k+1},p_{T,f,h}^{k+1})\in V_h\times Q_{T,h}$
be the discrete mechanical load lift defined by
\begin{align}
  \label{eq:tp:discrete-load-lift}
  a_u(u_{f,h}^{k+1},v_h)+b(p_{T,f,h}^{k+1},v_h)
    &=(\bar f^{k+1},v_h)_\Omega,\notag\\
  b(q_{T,h},u_{f,h}^{k+1})
    -\lambda^{-1}(p_{T,f,h}^{k+1},q_{T,h})_\Omega&=0
\end{align}
for all $(v_h,q_{T,h})\in V_h\times Q_{T,h}$.
The continuous lift $(\bar u_f^{k+1},\bar p_{T,f}^{k+1})\in V\times Q_T$
solves the same problem in the continuous spaces. By linearity of
\eqref{eq:tp:mechanical-load-lift}, it is the time average of the
load lift already defined in Section~\ref{sec:apriori}:
\[
  (\bar u_f^{k+1},\bar p_{T,f}^{k+1})
  =\frac{1}{\Delta t}\int_{I_{k+1}}(u_f(t),p_{T,f}(t))\,dt.
\]
Likewise, let $(u_{f,h}(t),p_{T,f,h}(t))$ solve
\eqref{eq:tp:discrete-load-lift} with load $f(t)$ in place of
$\bar f^{k+1}$; its time average is $(u_{f,h}^{k+1},p_{T,f,h}^{k+1})$.
Define the corresponding load functionals by
\begin{align}
  \label{eq:tp:load-functionals}
  G_f^{k+1}(q)&:=\lambda^{-1}\alpha(\bar p_{T,f}^{k+1},q)_\Omega,
    &&q\in Q,\notag\\
  G_{f,h}^{k+1}(q_h)&:=\lambda^{-1}\alpha(p_{T,f,h}^{k+1},q_h)_\Omega,
    &&q_h\in Q_h.
\end{align}

\begin{lemma}
  \label{lem:tp:scheme-equiv}
  The fully discrete scheme \eqref{eq:tp:FD-state} is equivalent,
  after eliminating $(u_h^{k+1},p_{T,h}^{k+1})$ via
  Lemma~\ref{lem:tp:discrete-mech}, to the reduced scalar scheme:
  find $p_h^{k+1}\in Q_h$ such that
  \begin{align}
    \label{eq:tp:reduced-scheme}
    m_{{\rm eff},h}(d_t p_h^{k+1},q_h)+(\kappa_p\nabla p_h^{k+1},\nabla q_h)_\Omega
    &=(\bar m_p^{k+1},q_h)_\Omega+G_{f,h}^{k+1}(q_h)
    \qquad\forall q_h\in Q_h,
  \end{align}
  where $\bar m_p^{k+1}:=\frac{1}{\Delta t}\int_{I_{k+1}}m_p(t)\,dt$.
\end{lemma}

\begin{proof}
  The first two equations of \eqref{eq:tp:FD-state} are the
  mechanical saddle-point problem for
  $(d_tu_h^{k+1},d_tp_{T,h}^{k+1})$, driven by $d_tp_h^{k+1}$
  and $\bar f^{k+1}$. By linearity and
  \eqref{eq:tp:discrete-load-lift},
  \begin{align}
    \label{eq:tp:mech-split}
    d_tu_h^{k+1}
      =\mathcal U_h(d_tp_h^{k+1})+u_{f,h}^{k+1},\qquad 
    d_tp_{T,h}^{k+1}
      =\mathcal P_{T,h}(d_tp_h^{k+1})+p_{T,f,h}^{k+1}.
  \end{align}
  Substitution into \eqref{eq:tp:FD-state3} and the definitions of
  $m_{{\rm eff},h}$ and $G_{f,h}^{k+1}$ give \eqref{eq:tp:reduced-scheme}.
  Conversely, given a solution of the reduced scheme and the initial
  mechanical variables, \eqref{eq:tp:mech-split} reconstructs the
  mechanical increments and hence the full discrete state.
  The load term depends only on $\bar f^{k+1}$. Its continuous counterpart
  is $G_f^{k+1}$, and their difference gives the load residual
  $R_F^{k+1}$ in the error equation below.
\end{proof}

\begin{lemma}
  \label{lem:tp:mass-consistency}
  Assume $\lambda\ge\lambda_0>0$ and the Stokes-type regularity
  assumption~(\ref{item:app:stokes}) of Appendix~\ref{app:regularity}.
  For $z\in H^1_0(\Omega)$ and $q_h\in Q_h$, 
  \begin{align}
    \label{eq:tp:mass-consistency}
    |m_{\rm eff}(z,q_h)-m_{{\rm eff},h}(P_hz,q_h)|\le Ch\|z\|_{H^1}\|q_h\|_0,
  \end{align}
  with $C$ independent of $h$, $\Delta t$, $s_0$, and uniform for
  $\lambda\ge\lambda_0$.
\end{lemma}

\begin{proof}
  By \eqref{eq:tp:reduced-mass-def}, \eqref{eq:tp:discrete-mass-def}
  and property~(a) of $P_h$, the storage part $s_1((I-P_h)z,q_h)_\Omega$
  vanishes, so
  \[
    m_{\rm eff}(z,q_h)-m_{{\rm eff},h}(P_hz,q_h)=-\lambda^{-1}\alpha\bigl((\mathcal P_T-\mathcal P_{T,h})z
      +\mathcal P_{T,h}(I-P_h)z,\,q_h\bigr)_\Omega .
  \]
  By Lemma~\ref{lem:tp:discrete-mech} and property~(b) of $P_h$,
  $\|\mathcal P_{T,h}(I-P_h)z\|_0\le C\|(I-P_h)z\|_0\le Ch\|z\|_{H^1}$.
  The pair $(\mathcal U_hz,\mathcal P_{T,h}z)$ is the Taylor--Hood
  Galerkin approximation of $(\mathcal Uz,\mathcal P_Tz)$, so
  Babu\v{s}ka--Brezzi quasi-optimality, the Taylor--Hood approximation
  properties (Definition~\ref{def:tp:interp}), and assumption~(\ref{item:app:stokes}) with $F=0$,
  $\zeta=-\lambda^{-1}\alpha z\in H^1_0(\Omega)$ give
  \[
    \|(\mathcal P_T-\mathcal P_{T,h})z\|_0
    \le Ch\bigl(\|\mathcal Uz\|_{H^2}+\|\mathcal P_Tz\|_{H^1}\bigr)
    \le Ch\|z\|_{H^1}.
  \]
  Using $\lambda^{-1}\le\lambda_0^{-1}$ completes the proof.
\end{proof}

\begin{lemma}[Stability of the reduced backward dual]
  \label{lem:tp:reduced-dual-stability}
  Let $\varphi^n\in Q_h$, $n=1,\ldots,N$, with $\varphi^{N+1}:=0$,
  solve
  \begin{align}
    \label{eq:tp:dual-reduced}
    m_{{\rm eff},h}\!\left(\frac{\varphi^n-\varphi^{n+1}}{\Delta t},q_h\right)
    +(\kappa_p\nabla\varphi^n,\nabla q_h)_\Omega
    =(\xi^n,q_h)_{L^2}\quad\forall q_h\in Q_h,
  \end{align}
  for a given sequence $\xi^n\in Q_h$.  Then
  \begin{align}
    \label{eq:tp:dual-stab-L2}
    \Bigl(\Delta t\sum_{n=1}^N\|\varphi^n\|_0^2\Bigr)^{1/2}
    &\le C\Bigl(\Delta t\sum_{n=1}^N\|\xi^n\|_0^2\Bigr)^{1/2},
    \\
    \label{eq:tp:dual-stab-H1}
    \Bigl(\Delta t\sum_{n=1}^N\|\nabla\varphi^n\|_0^2\Bigr)^{1/2}
    &\le C\Bigl(\Delta t\sum_{n=1}^N\|\xi^n\|_0^2\Bigr)^{1/2},
    \\
    \label{eq:tp:dual-stab-dt}
    \Bigl(\Delta t\sum_{n=1}^N
    \left\|\frac{\varphi^n-\varphi^{n+1}}{\Delta t}\right\|_0^2
    \Bigr)^{1/2}
    &\le C\Bigl(\Delta t\sum_{n=1}^N\|\xi^n\|_0^2\Bigr)^{1/2},
    \\
    \label{eq:tp:dual-stab-max}
    \max_{1\le n\le N}\|\varphi^n\|_0
    &\le C\Bigl(\Delta t\sum_{n=1}^N\|\xi^n\|_0^2\Bigr)^{1/2},
  \end{align}
  with $C$ independent of $h$ and $\Delta t$.  It depends on
  the coercivity of $m_{{\rm eff},h}$ and may therefore depend on $\lambda$
  when $s_0=0$.  For $s_0\ge s_*>0$, the estimates are
  $\lambda$-uniform under the stated diffusion assumptions.
\end{lemma}

\begin{proof}
  Set $\delta_{\tau}\varphi^n := (\varphi^n - \varphi^{n+1})/\Delta t$.
  For \eqref{eq:tp:dual-stab-L2} and \eqref{eq:tp:dual-stab-H1} test \eqref{eq:tp:dual-reduced} with $q_h=\varphi^n$.
  Since $m_{{\rm eff},h}$ is symmetric and coercive on $Q_h$, it induces the
  norm $\|\varphi\|_{m_{{\rm eff},h}}:=m_{{\rm eff},h}(\varphi,\varphi)^{1/2}$, and
  \[
    m_{{\rm eff},h}(\delta_{\tau}\varphi^n,\varphi^n)
    =\frac{1}{2\Delta t}\bigl(
      \|\varphi^n\|_{m_{{\rm eff},h}}^2
      -\|\varphi^{n+1}\|_{m_{{\rm eff},h}}^2
      +\|\varphi^n-\varphi^{n+1}\|_{m_{{\rm eff},h}}^2
    \bigr)
    \ge\frac{1}{2\Delta t}\bigl(
      \|\varphi^n\|_{m_{{\rm eff},h}}^2-\|\varphi^{n+1}\|_{m_{{\rm eff},h}}^2
    \bigr).
  \]
  Using $(\kappa_p\nabla\varphi^n,\nabla\varphi^n)_\Omega=\|\varphi^n\|_Q^2$,
  Poincar\'e's inequality on $Q_h$, and
  Young's inequality, we have
  \[
    (\xi^n,\varphi^n)\le C\|\xi^n\|_0\|\varphi^n\|_{H^1}
    \le\varepsilon\|\varphi^n\|_{H^1}^2
       +C_\varepsilon\|\xi^n\|_0^2.
  \]
  Multiplying by $\Delta t$, summing $n=1,\ldots,N$, using
  $\varphi^{N+1}=0$, and choosing $\varepsilon$ small:
  \begin{align*}
    \frac{1}{2}\|\varphi^1\|_{m_{{\rm eff},h}}^2
    +c\Delta t\sum_{n=1}^N\|\varphi^n\|_{H^1}^2
    \le C\Delta t\sum_{n=1}^N\|\xi^n\|_0^2.
  \end{align*}
  The $m_{{\rm eff},h}$-coercivity and Poincar\'e inequality give
  \eqref{eq:tp:dual-stab-L2}--\eqref{eq:tp:dual-stab-H1}.

  For \eqref{eq:tp:dual-stab-max}, for arbitrary $1\le m\le N$, sum the energy inequality from $n=m$ to $N$.  Using $\varphi^{N+1}=0$:
  \[
    \frac{1}{2}\|\varphi^m\|_{m_{{\rm eff},h}}^2
    +c\Delta t\sum_{n=m}^N\|\varphi^n\|_{H^1}^2
    \le C\Delta t\sum_{n=m}^N\|\xi^n\|_0^2
    \le C\Delta t\sum_{n=1}^N\|\xi^n\|_0^2.
  \]
  Taking the maximum over $m=1,\ldots,N$ and using $m_{{\rm eff},h}$-coercivity:
  \begin{align*}
    \max_{1\le n\le N}\|\varphi^n\|_0^2
    \le C\max_{1\le m\le N}\|\varphi^m\|_{m_{{\rm eff},h}}^2
    \le C\Delta t\sum_{n=1}^N\|\xi^n\|_0^2,
  \end{align*}
  giving \eqref{eq:tp:dual-stab-max}.

  For \eqref{eq:tp:dual-stab-dt}, test \eqref{eq:tp:dual-reduced} with $q_h=\delta_{\tau}\varphi^n$.
  By symmetry of $(\kappa_p\nabla\cdot,\nabla\cdot)_\Omega$:
  \[
    (\kappa_p\nabla\varphi^n,\nabla \delta_{\tau}\varphi^n)_\Omega
    =\frac{1}{2\Delta t}\bigl(
      \|\varphi^n\|_Q^2
      -\|\varphi^{n+1}\|_Q^2
      +\|\varphi^n-\varphi^{n+1}\|_Q^2
    \bigr)
    \ge\frac{1}{2\Delta t}\bigl(
      \|\varphi^n\|_Q^2-\|\varphi^{n+1}\|_Q^2
    \bigr).
  \]
  Multiplying by $\Delta t$, summing $n=1,\ldots,N$, using
  $\varphi^{N+1}=0$, $m_{{\rm eff},h}$-coercivity, and Young's inequality:
  \[
    c\Delta t\sum_{n=1}^N\|\delta_{\tau}\varphi^n\|_0^2
    \le\Delta t\sum_{n=1}^N\|\xi^n\|_0\|\delta_{\tau}\varphi^n\|_0
    \le\frac{c}{2}\Delta t\sum_{n=1}^N\|\delta_{\tau}\varphi^n\|_0^2
    +C\Delta t\sum_{n=1}^N\|\xi^n\|_0^2,
  \]
  giving \eqref{eq:tp:dual-stab-dt}.
\end{proof}

\begin{lemma}
  \label{lem:tp:reduced-error}
  Under the regularity \eqref{eq:tp:reg-state-adj} and the initial
  approximation assumption $\|p_h^0-p(0)\|_0=O(h)$, the discrete
  solution $p_h^j$ from \eqref{eq:tp:reduced-scheme}
  satisfies
  \begin{align}
    \label{eq:tp:reduced-error-bound}
    \Bigl(\sum_{j=1}^N\int_{I_j}\|p(t)-p_h^j\|_0^2\,dt\Bigr)^{1/2}
    \le C(h+\Delta t),
  \end{align}
  with $C$ independent of $h$ and $\Delta t$ but depending on the lower bounds of $c_{\rm eff}$ and $c_{{\rm eff},h}$.
\end{lemma}

\begin{proof}
  Set $\bar p^0:=p(0)$.  For $k=1,\ldots,N$, let
  $\bar p^k:=\frac{1}{\Delta t}\int_{I_k}p(t)\,dt$ and
  $d^k:=p(t_k)-\bar p^k$, so $d^0=0$.

  For $t\in I_j$ we decompose $p(t)-p_h^j$ into 
  \[
    p(t)-p_h^j
    = (I-P_h)p(t)
    + P_h(p(t)-\bar p^j)
    + (P_h\bar p^j-p_h^j) =: \eta(t) + \mu^j(t) + \xi^j.
  \]
  By property~(b) in Lemma~\ref{lem:interpolation-property},  $\int_0^T\|\eta(t)\|_0^2\,dt\le Ch^2\|p\|_{L^2_tH^1_x}^2$. We can also obtain $\sum_{j=1}^N\int_{I_j}\|\mu^j\|_0^2\,dt\le
  C\Delta t^2\|\dot p\|_{L^2_tL^2_x}^2$ because 
  $\|p(t)-\bar p^j\|_0\le\int_{I_j}\|\dot p\|_0\,dt$. Thus, $\|\eta \|_{L_t^2 L_x^2}^2 + \sum_{j=1}^N\int_{I_j}\|\mu^j\|_0^2\,dt \le C(h^2+\Delta t^2)$.

  It remains to estimate $\Delta t \sum_{n=1}^N \|\xi^n\|_{0}^2$. For this we derive an error equation and apply discrete duality. 
  
  \emph{Error equation of $\xi^j$}:
  For the error equation for $\xi^k=P_h\bar p^k-p_h^k\in Q_h$, 
  since $d_t(P_h\bar p^{k+1})=P_hd_t\bar p^{k+1}$ and $d_t\bar p^{k+1}
  =\frac{\bar p^{k+1}-\bar p^k}{\Delta t}$, subtracting
  \eqref{eq:tp:reduced-scheme} from the time-averaged continuous
  equation gives, for all $q_h\in Q_h$:
  \begin{align}
    \label{eq:tp:xi-error-dg}
    m_{{\rm eff},h}(d_t\xi^{k+1},q_h)+(\kappa_p\nabla\xi^{k+1},\nabla q_h)_\Omega
    =-R_M^{k+1}(q_h)-R_S^{k+1}(q_h)-C_\tau^{k+1}(q_h)+R_F^{k+1}(q_h),
  \end{align}
  where
  \begin{align*}
    R_M^{k+1}(q_h)&:=m_{\rm eff}(\bar{\dot p}^{k+1},q_h)
      -m_{{\rm eff},h}(P_h\bar{\dot p}^{k+1},q_h),
    \\
    R_S^{k+1}(q_h)&:=\frac{1}{\Delta t}\int_{I_{k+1}}(\kappa_p\nabla p(t),\nabla q_h)_\Omega\,dt
      -(\kappa_p\nabla P_h\bar p^{k+1},\nabla q_h)_\Omega,
    \\
    C_\tau^{k+1}(q_h)&:=m_{{\rm eff},h}\!\left(P_h\!\left(\frac{d^{k+1}-d^k}{\Delta t}\right),q_h\right),
    \\
    R_F^{k+1}(q_h)&:=\lambda^{-1}\alpha
      \bigl(\bar p_{T,f}^{k+1}-p_{T,f,h}^{k+1},\,q_h\bigr)_\Omega,
  \end{align*}
  and $\bar{\dot p}^{k+1}:=\frac{1}{\Delta t}\int_{I_{k+1}}\dot p(t)\,dt
  =\frac{p(t_{k+1})-p(t_k)}{\Delta t}$.
  The load lifts in \eqref{eq:tp:discrete-load-lift} and its
  continuous counterpart give
  $R_F^{k+1}=G_f^{k+1}|_{Q_h}-G_{f,h}^{k+1}$ by
  \eqref{eq:tp:load-functionals}. This residual is the mechanical
  Galerkin error in the total pressure load lift.
  Set $\xi^0:=P_hp(0)-p_h^0$.
  Since $p_h^0\in Q_h$ we have $P_hp_h^0=p_h^0$, so
  $\xi^0=P_h(p(0)-p_h^0)$ and hence $\|\xi^0\|_0\le\|p(0)-p_h^0\|_0=O(h)$
  by the initial approximation assumption.

  \emph{An equation of $\Delta t \sum_{n=1}^N \|\xi^n\|_{0}^2$ via the dual equation}:
  Consider a discrete dual problem:
  Let $\varphi^n\in Q_h$, $n=1,\ldots,N$, with $\varphi^{N+1}:=0$,
  solve \eqref{eq:tp:dual-reduced} with right-hand side $\xi^n$. Then,
  Lemma~\ref{lem:tp:reduced-dual-stability} gives
  \eqref{eq:tp:dual-stab-L2}--\eqref{eq:tp:dual-stab-max}.
  Testing \eqref{eq:tp:xi-error-dg} at step $n$ with $q_h=\varphi^n$, summing $\Delta t \sum_{n=1}^N$, and applying summation by parts using $\varphi^{N+1}=0$, 
  \begin{align}
    \label{eq:tp:duality-xi}
    \Delta t\sum_{n=1}^N\|\xi^n\|_0^2
    =-\Delta t\sum_{n=1}^{N}\bigl[
      R_M^{n}(\varphi^{n})
      +R_S^{n}(\varphi^{n})
      +C_\tau^{n}(\varphi^{n})
      -R_F^{n}(\varphi^{n})\bigr]
    +m_{{\rm eff},h}(\xi^0,\varphi^1).
  \end{align}
  The boundary term $m_{{\rm eff},h}(\xi^0,\varphi^1)$ can be estimated by 
  \begin{align*}
  |m_{{\rm eff},h}(\xi^0,\varphi^1)|\le C\|\xi^0\|_0\|\varphi^1\|_0
  \le C\|\xi^0\|_0\max_{1\le n\le N}\|\varphi^n\|_0 \le  Ch\left(\Delta t\sum_{n=1}^N\|\xi^n\|_0^2\right)^{1/2}     
  \end{align*}
  from the previous estimate of $\|\xi^0\|_0$ and \eqref{eq:tp:dual-stab-max}.

  \emph{Estimates of residuals}: 
  We now prove error bounds of the residual terms $R_M$, $R_S$, $C_{\tau}$, $R_F$. 

  We first estimate $R_S$. Since the diffusion form is time-independent and linear,
  and $\bar p^{k+1}$ is the exact cell average of $p$ over $I_{k+1}$:
  \[
    \frac{1}{\Delta t}\int_{I_{k+1}}(\kappa_p\nabla p(t),\nabla \varphi^{k+1})_\Omega\,dt
    = (\kappa_p\nabla\bar p^{k+1},\nabla \varphi^{k+1})_\Omega.
  \]
  Hence $R_S^{k+1}(\varphi^{k+1})
  = (\kappa_p\nabla(I-P_h)\bar p^{k+1},\nabla \varphi^{k+1})_\Omega$.
  Since $\bar p^{k+1}\in H^2(\Omega)$ by Jensen's inequality and $p\in L^2_tH^2_x$,
  property~(c) in Lemma~\ref{lem:interpolation-property} gives $\|(I-P_h)\bar p^{k+1}\|_{H^1}\le Ch|\bar p^{k+1}|_{H^2}$.
  Multiplying by $\Delta t$, summing, Cauchy--Schwarz and H\"{o}lder inequalities with \eqref{eq:tp:dual-stab-H1} give
  \begin{align}
    \label{eq:R_S-estimate}
    \left| \Delta t\sum_{n=1}^N R_S^{n}(\varphi^{n}) \right|\le Ch\|p\|_{L^2_tH^2_x} \left( \Delta t \sum_{n=1}^N \|\xi^n\|_0^2 \right)^{\frac 12}.
  \end{align}
  We now estimate $R_F$. The continuous load-driven system for $(\bar u_f^{k+1},\bar p_{T,f}^{k+1})$
  and its discrete counterpart $(u_{f,h}^{k+1},p_{T,f,h}^{k+1})$ both solve
  a perturbed Stokes problem with the same right-hand side $(\bar f^{k+1},0)$
  as in \eqref{eq:tp:discrete-load-lift}. By quasi-optimality for the Taylor--Hood pair, which bounds the pressure
  error by the best-approximation errors of both the displacement and
  the total pressure (in Definition~\ref{def:tp:interp}), and the Stokes-type regularity
  assumption~(\ref{item:app:stokes}) of Appendix~\ref{app:regularity}
  (with $F=\bar f^{k+1}$, $\zeta=0$),
  \[
    \|\bar p_{T,f}^{k+1}-p_{T,f,h}^{k+1}\|_0
    \le Ch\bigl(\|\bar u_f^{k+1}\|_{H^2}+\|\bar p_{T,f}^{k+1}\|_{H^1}\bigr)
    \le Ch\|\bar f^{k+1}\|_0.
  \]
  Hence, for the dual solution $\{\varphi^n\} \in Q_h$,
  \begin{align}
    \label{eq:tp:R_F-estimate}
    \left|\Delta t\sum_{n=1}^N R_F^n(\varphi^n)\right|
    \le C\lambda^{-1}\alpha\,h\,\|f\|_{L^2_tL^2_x}
    \Bigl(\Delta t\sum_{n=1}^N\|\varphi^n\|_0^2\Bigr)^{1/2}
    \le Ch\cdot C_{\rm data}\cdot  \left( \Delta t \sum_{n=1}^N \|\xi^n\|_0^2 \right)^{\frac 12} 
  \end{align}
  by \eqref{eq:tp:dual-stab-L2} 
  where $C_{\rm data}:=\lambda^{-1}\alpha T^{1/2}\|f\|_{L^2_tL^2_x}$
  is $\lambda$-uniform for $\lambda\ge\lambda_0$.

  


  We prove bounds for $R_M$ and $C_\tau$. Let $E_h(z,q_h):= m_{\rm eff}(z,q_h)-m_{{\rm eff},h}(P_h z,q_h)$ for $z \in H_0^1(\Omega)$.
  Since $\bar{\dot p}^{\,n}=d_t\bar p^n+d_td^n$ with
  $d_td^n:=(d^n-d^{n-1})/\Delta t$, the definitions give
  \[
    R_M^n(q_h)+C_\tau^n(q_h)=E_h(d_t\bar p^n,q_h) 
    +m_{\rm eff}(d_td^n,q_h).
  \]
  Summation by parts with $\delta_{\tau}\varphi^n=(\varphi^n-\varphi^{n+1})/\Delta t$,
  $\varphi^{N+1}=0$, $\bar p^0=p(0)$, and $d^0=0$ gives
  \[
    \Delta t\sum_{n=1}^N(R_M^n+C_\tau^n)(\varphi^n)
    =\Delta t\sum_{n=1}^N\bigl[E_h(\bar p^n,\delta_{\tau}\varphi^n)
      +m_{\rm eff}(d^n,\delta_{\tau}\varphi^n)\bigr]-E_h(p(0),\varphi^1).
  \]
  Here $\bar p^n\in H^1_0(\Omega)$, and $p(0)\in H^1_0(\Omega)$ by the
  embedding $H^1(0,T;L^2)\cap L^2(0,T;H^2\cap H^1_0)\subset
  C([0,T];H^1_0)$.  Writing
  $d^n=\Delta t^{-1}\int_{I_n}\int_t^{t_n}\dot p(s)\,ds\,dt$,
  Cauchy--Schwarz yields
  $\Delta t\sum_{n=1}^N\|d^n\|_0^2\le\Delta t^2\|\dot p\|_{L^2_tL^2_x}^2$.
  By Lemma~\ref{lem:tp:mass-consistency},
  $\Delta t\sum_n\|\bar p^n\|_{H^1}^2\le\|p\|_{L^2_tH^1_x}^2$ due to Jensen's inequality,
  boundedness of $m_{\rm eff}$,
  \eqref{eq:tp:dual-stab-dt} and \eqref{eq:tp:dual-stab-max} yield
  \begin{align}
    \label{eq:R_M-C_tau-estimate}
    \Bigl|\Delta t\sum_{n=1}^N(R_M^n+C_\tau^n)(\varphi^n)\Bigr|
    \le C\bigl(h\|p\|_{L^2_tH^1_x}+h\|p(0)\|_{H^1}
      +\Delta t\|\dot p\|_{L^2_tL^2_x}\bigr)
      \Bigl(\Delta t\sum_{n=1}^N\|\xi^n\|_0^2\Bigr)^{1/2}.
  \end{align}
  The conclusion then follows by collecting \eqref{eq:R_S-estimate}, \eqref{eq:tp:R_F-estimate}, \eqref{eq:R_M-C_tau-estimate}, and the triangle inequality.
\end{proof}

\begin{lemma}
  \label{lem:tp:mech-recovery}
  Under the assumptions of Lemma~\ref{lem:tp:reduced-error},
  \eqref{eq:tp:fixed-state-mech} holds.
\end{lemma}

\begin{proof}
  We use the integrated mechanical reconstruction.  For
  $t\in[0,T]$ and $j=1,\ldots,N$,
  \begin{align*}
    u(t)&=u(0)+\int_0^t u_f(s)\,ds
             +\mathcal U(p(t)-p(0)),\\
    p_T(t)&=p_T(0)+\int_0^t p_{T,f}(s)\,ds
             +\mathcal P_T(p(t)-p(0)),\\
    u_h^j&=u_h^0+\Delta t\sum_{\ell=1}^j u_{f,h}^{\ell}
             +\mathcal U_h(p_h^j-p_h^0),\\
    p_{T,h}^j&=p_{T,h}^0
             +\Delta t\sum_{\ell=1}^j p_{T,f,h}^{\ell}
             +\mathcal P_{T,h}(p_h^j-p_h^0).
  \end{align*}
  Here $(u_{f,h}^{\ell},p_{T,f,h}^{\ell})$ is the mechanical
  load lift in \eqref{eq:tp:discrete-load-lift}. Taking the cell average
  of the first two identities and comparing with the last two
  reduces the error to three contributions.

  First, the initial-data errors are $O(h)$ by the numerical
  initial approximation assumption.  Second, the continuous and
  discrete load lifts satisfy, by mixed Galerkin quasi-optimality
  and the mechanical $H^2\times H^1$ regularity estimate,
  \[
    \|u_f(s)-u_{f,h}(s)\|_V
       +\|p_{T,f}(s)-p_{T,f,h}(s)\|_0
       \le Ch\|f(s)\|_0,
  \]
  The time averages of these lifts are
  $u_{f,h}^{\ell},p_{T,f,h}^{\ell}$, as defined above.
  Replacing a cumulative integral at $t_j$ by its average over $I_j$
  produces an additional $O(\Delta t)$ term in $\ell^2_t$, since
  \[
    \left\|\frac1{\Delta t}\int_{I_j}\int_t^{t_j}u_f(s)\,ds\,dt\right\|_V
       \le\Delta t^{1/2}\|u_f\|_{L^2(I_j;V)},
  \]
  and similarly for $p_{T,f}$.  Consequently the total load-lift
  contribution is bounded in $\ell^2_t(V\times Q_T)$ by
  $C(h+\Delta t)\|f\|_{L^2_tL^2_x}$.

  Third, for $y\in H^1_0(\Omega)$, mixed Galerkin approximation,
  the $L^2$-approximation property of $P_h$, and the regularity
  assumptions give
  \begin{align}
    \label{eq:tp:mech-operator-approx}
    \|\mathcal Uy-\mathcal U_hP_hy\|_V
       +\|\mathcal P_Ty-\mathcal P_{T,h}P_hy\|_0
       \le Ch\|y\|_{H^1}.
  \end{align}
  Apply this to $y=\bar p^j-p(0)$.  By the $L^2_tH^2_x\cap
  H^1_tL^2_x$ regularity in Assumption~R, the embedding $L_t^2 H_x^2 \cap H_t^1 L_x^2 \hookrightarrow C^0([0,T];H^1(\Omega))$ gives
  $\|p\|_{C^0([0,T];H^1)}<\infty$.
  Using the uniform stability of $(\mathcal U_h,\mathcal P_{T,h})$,
  \begin{align*}
    &\|\mathcal U(\bar p^j-p(0))
           -\mathcal U_h(p_h^j-p_h^0)\|_V
      +\|\mathcal P_T(\bar p^j-p(0))
           -\mathcal P_{T,h}(p_h^j-p_h^0)\|_0\\
    &\qquad\le Ch\|\bar p^j-p(0)\|_{H^1}
       +C\bigl(\|P_h\bar p^j-p_h^j\|_0
                   +\|P_hp(0)-p_h^0\|_0\bigr).
  \end{align*}
  The first and last terms are $O(h)$ in $\ell^2_t$.
  For the middle term, the decomposition in the proof of
  Lemma~\ref{lem:tp:reduced-error} and Jensen's inequality give
  \[
    \Delta t\sum_{j=1}^N\|P_h\bar p^j-p_h^j\|_0^2
       \le \sum_{j=1}^N\int_{I_j}\|p(t)-p_h^j\|_0^2\,dt
       \le C(h+\Delta t)^2.
  \]
  Combining the three contributions proves
  \eqref{eq:tp:fixed-state-mech}.
\end{proof}

\begin{proof}[Proof of Lemma~\ref{lem:tp:fixed-control-state-error}]
  Estimate \eqref{eq:tp:fixed-state-diff} follows directly from
  Lemma~\ref{lem:tp:reduced-error}.
  The reduction to the scalar scheme \eqref{eq:tp:reduced-scheme} is
  established in Lemmas~\ref{lem:tp:mech-operator},
  \ref{lem:tp:reduced-mass}, \ref{lem:tp:discrete-mech},
  and~\ref{lem:tp:scheme-equiv}.
  Estimate \eqref{eq:tp:fixed-state-mech} follows from
  Lemma~\ref{lem:tp:mech-recovery}.
\end{proof}

\subsection{Estimates for the control error}
\label{sec:THM:control-error}

We first record a strong monotonicity property of the discrete
reduced gradient $j_h'$ on $\mathcal M_{p,ad}^h$, and then a fixed-control estimate for the
adjoint error, before combining them to bound the control error.

\begin{lemma}
  \label{lem:tp:Jh-monotone}
  For any $m_{p,h},\tilde m_{p,h}\in\mathcal{M}_{p,ad}^h$,
  \begin{align}
    \label{eq:tp:Jh-monotone}
    \bigl(j_h'(m_{p,h}) - j_h'(\tilde m_{p,h}),\,
      m_{p,h}-\tilde m_{p,h}\bigr)_{L^2_tL^2_x}
    \ge \gamma_p \|m_{p,h}-\tilde m_{p,h}\|_{L^2_tL^2_x}^2,
  \end{align}
  where $j_h'$ is the discrete gradient \eqref{eq:tp:Jh-gradient}.
\end{lemma}

\begin{proof}
  Let $m_{p,d} := m_{p,h} - \tilde m_{p,h}$ and let $(u_{d,h}^{k+1}, p_{T,d,h}^{k+1},
  p_{d,h}^{k+1})$ denote the corresponding
  difference of the discrete states
  $u_h^{k+1}(m_{p,h}) - u_h^{k+1}(\tilde m_{p,h})$, etc.
  The derivative of the affine control-to-state map in the direction
  $m_{p,d}$ is $(u_{d,h},p_{T,d,h},p_{d,h})$, independent of the base
  point.  Since $j_h$ is quadratic, we obtain 
  \begin{align*}
    \bigl(j_h'(m_{p,h}) - j_h'(\tilde m_{p,h}),\,
      m_{p,d}\bigr)_{L^2_tL^2_x} =\alpha_{u,C}\Delta t \sum_{k=0}^{N-1}\|u_{d,h}^{k+1}\|_0^2
      +\alpha_{p,C}\Delta t \sum_{k=0}^{N-1}\|p_{d,h}^{k+1}\|_0^2
      +\gamma_p \|m_{p,d}\|_{L^2_tL^2_x}^2.
  \end{align*}
  Since $\alpha_{u,C},\alpha_{p,C} \ge 0$ and $\gamma_p > 0$,
  dropping the nonnegative tracking terms yields
  \eqref{eq:tp:Jh-monotone}.
\end{proof}

We next bound the adjoint error for a \emph{fixed} control, which is
the key ingredient for the control-error estimate.

\begin{lemma}
  \label{lem:tp:adjoint-error-fixed-control}
  Let $\widehat m_p \in \mathcal{M}_{p,ad}$ be the
  continuous optimal control, with continuous state
  $(\widehat u, \widehat p_T, \widehat p)$ and continuous adjoint
  $(\widehat v^*, \widehat q_T^*, \widehat q^{\,*})$ solving
  \eqref{eq:tp:adjoint}--\eqref{eq:tp:adj:TC}.  Let
  $(u_h(\widehat m_p), p_{T,h}(\widehat m_p), p_h(\widehat m_p))$ be
  the auxiliary fully discrete state driven by $\widehat m_p$
  as in \eqref{eq:tp:FD-state}, and let
  $(v_h^*(\widehat m_p), q_{T,h}^*(\widehat m_p), q_h^*(\widehat m_p))$ be
  the corresponding auxiliary fully discrete adjoint, defined by
  \eqref{eq:tp:FD-adj}--\eqref{eq:tp:FD-adj-TC} with the right-hand
  sides of \eqref{eq:tp:FD-adj1} and \eqref{eq:tp:FD-adj3}
  and with homogeneous terminal condition
  \eqref{eq:tp:FD-adj-TC}.
  Under  \eqref{eq:tp:reg-state-adj} in Assumption~R, there holds
  \begin{align}
    \label{eq:tp:adjoint-error-fixed-control}
    \Bigl(\Delta t\sum_{k=0}^{N-1}
      \|q_h^{*,k}(\widehat m_p) - \overline{\widehat q^{\,*}}{}^{\,k+1}\|_0^2
    \Bigr)^{1/2}
    \le C(h+\Delta t),
  \end{align}
  where $\overline{\widehat q^{\,*}}{}^{\,k+1}$ denotes the time-average
  of $\widehat q^{\,*}$ over $(t_k,t_{k+1})$, and $C>0$ is
  independent of $h$ and $\Delta t$, uniform for $\lambda\ge\lambda_0>0$ if $s_0 \ge s_* >0$ and the assumed regularity bounds are uniform.
\end{lemma}

\begin{proof}
  Introduce the auxiliary discrete adjoint $(\widetilde v_h^{*,k},
  \widetilde q_{T,h}^{*,k}, \widetilde q_h^{*,k})$ defined by \eqref{eq:tp:FD-adj}--\eqref{eq:tp:FD-adj-TC}
  with the interior right-hand sides of \eqref{eq:tp:FD-adj1} and
  \eqref{eq:tp:FD-adj3} driven by the
  \emph{continuous} optimal state $(\widehat u, \widehat p)$
  (i.e.\ $\alpha_{u,C}(\overline{\widehat u}{}^{\,k+1} - \bar u_C^{k+1},v_h)_\Omega$ in
  place of $\alpha_{u,C}(u_h^{k+1}(\widehat m_p) -
  \bar u_C^{k+1},v_h)_\Omega$, and analogously for the $p$-equation).
  Both auxiliary adjoint systems use the same
  homogeneous terminal condition \eqref{eq:tp:FD-adj-TC}; hence their
  difference has homogeneous terminal data.  We split
  \begin{align*}
    q_h^{*,k}(\widehat m_p) - \overline{\widehat q^{\,*}}{}^{\,k+1}
    = \bigl(q_h^{*,k}(\widehat m_p) - \widetilde q_h^{*,k}\bigr)
    + \bigl(\widetilde q_h^{*,k} - \overline{\widehat q^{\,*}}{}^{\,k+1}\bigr) =: E_1^{q,k} + E_2^{q,k}.
  \end{align*}

  As the first step, we prove a bound on $E_1$.
  Subtracting the two discrete adjoint systems (driven by
  $u_h^{k+1}(\widehat m_p)$ versus $\overline{\widehat u}{}^{\,k+1}$, etc.\ in
  the interior right-hand sides, with homogeneous terminal
  data) shows that
  $(E_1^{v},E_1^{q_T},E_1^{q}) :=
  (v_h^*(\widehat m_p) - \widetilde v_h^*,\,
   q_{T,h}^*(\widehat m_p) - \widetilde q_{T,h}^*,\,
   q_h^*(\widehat m_p) - \widetilde q_h^*)$
  solves the discrete adjoint system \eqref{eq:tp:FD-adj} with
  interior right-hand sides
  \begin{align*}
    F_u^{k+1} = u_h^{k+1}(\widehat m_p) - \overline{\widehat u}{}^{\,k+1},
    \qquad
    F_p^{k+1} = p_h^{k+1}(\widehat m_p) - \overline{\widehat p}{}^{\,k+1},
  \end{align*}
  and with homogeneous terminal condition \eqref{eq:tp:FD-adj-TC}.
  Lemma~\ref{lem:tp:FD-adj-stability} directly gives
  \begin{align}
    \label{eq:tp:E1-bound}
    \Bigl(\Delta t\sum_{k=0}^{N-1}\|E_1^{q,k}\|_0^2\Bigr)^{1/2}
    &\le C \Bigl(\Delta t\sum_{j=1}^{N}\bigl(
        \|u_h^j(\widehat m_p) - \overline{\widehat u}{}^{\,j}\|_0^2
        + \|p_h^j(\widehat m_p) - \overline{\widehat p}{}^{\,j}\|_0^2
      \bigr)\Bigr)^{1/2}.
  \end{align}
  By Lemma~\ref{lem:tp:fixed-control-state-error}
  (applied to the auxiliary state driven by $\widehat m_p$),
  the right-hand side is bounded by $C(h+\Delta t)$.
  Therefore
  \begin{align}
    \label{eq:tp:E1-final}
    \Bigl(\Delta t\sum_{k=0}^{N-1}\|E_1^{q,k}\|_0^2\Bigr)^{1/2}
    \le C(h+\Delta t).
  \end{align}

  We now prove a bound on $E_2$.
  Reverse time as in the proof of Lemma~\ref{lem:tp:FD-adj-stability}:
  $\check q_h^l:=\widetilde q_h^{*,N-l}$ and
  $\check q^{\,*}(\tau):=\widehat q^{\,*}(T-\tau)$.  The interval
  $(t_k,t_{k+1}]$ corresponds to the reversed cell $I_l$ with $l=N-k$
  (up to endpoints), so
  $\overline{\widehat q^{\,*}}{}^{\,k+1}=\overline{\check q^{\,*}}{}^{\,l}$
  and $E_2^{q,k}=\check q_h^l-\overline{\check q^{\,*}}{}^{\,l}$.
  By \eqref{eq:tp:adj-reversed}, eliminating the mechanical variables
  yields a scheme of the form \eqref{eq:tp:reduced-scheme} for
  $\check q_h^l$, with the same effective mass $m_{{\rm eff},h}$, zero
  initial value, and data given by the cell averages of the adjoint
  forcing evaluated at the continuous optimal state; the continuous
  reversed adjoint satisfies the corresponding reduced equation with
  $m_{\rm eff}$ and zero initial value.  The load-lift term enters with
  the opposite sign, which does not affect the bound on $R_F$.  Under
  Assumption~R, Lemma~\ref{lem:tp:reduced-error} therefore applies to
  $(\check q^{\,*},\check q_h)$, and Jensen's inequality,
  $\Delta t\|\overline{\check q^{\,*}}{}^{\,l}-\check q_h^l\|_0^2
  \le\int_{I_l}\|\check q^{\,*}(\tau)-\check q_h^l\|_0^2\,d\tau$, gives
  \begin{align}
    \label{eq:tp:E2-final}
    \Bigl(\Delta t\sum_{k=0}^{N-1}\|E_2^{q,k}\|_0^2\Bigr)^{1/2}
    \le C(h+\Delta t).
  \end{align}

  Finally, by the triangle inequality,
  \eqref{eq:tp:adjoint-error-fixed-control} follows from
  \eqref{eq:tp:E1-final} and \eqref{eq:tp:E2-final}.
\end{proof}

\subsection{Main error estimate for the control}
\label{sec:THM:main-error}

\begin{theorem}
  \label{thm:tp:control-error}
  Let $\widehat{m}_p$ and $\widehat{m}_{p,h}$
  be the continuous and discrete optimal controls, respectively.
  Assume $s_0\ge s_*>0$.  Under Assumption~R
  \eqref{eq:tp:reg-state-adj}, with a regularity
  bound $C_{\rm reg}$ uniform for $\lambda\ge\lambda_0>0$ (depending on
  $\lambda_0$), and with independence of $C_{\rm reg}$ from $\kappa_p$
  requiring the elliptic regularity constants in
  \eqref{eq:tp:reg-state-adj} to be uniform in $\kappa_p$,
  there holds
  \begin{align}
    \label{eq:tp:control-error}
    \|\widehat{m}_p - \widehat{m}_{p,h}\|_{L^2_tL^2_x}
    \le C(h+\Delta t),
  \end{align}
  where $C > 0$ depends on $\gamma_p,T$ and $C_{\rm
  reg}$, is independent of $h$ and $\Delta t$, and is uniform for
  $\lambda\ge\lambda_0>0$.
\end{theorem}

\begin{rmk}
  \label{rmk:tp:control-regularity}
  Assumption~R implies $\widehat q^{\,*}\in H^1((0,T)\times\Omega)$.
  For the constant bounds $l_p,L_p$, the scalar box projection is
  $1$-Lipschitz, so the projection formula \eqref{eq:tp:proj-p} and
  the Sobolev chain rule therefore give
  \begin{align}
    \label{eq:tp:control-regularity}
    \widehat m_p\in H^1(0,T;L^2(\Omega))\cap L^2(0,T;H^1(\Omega)).
  \end{align}
  More precisely, its weak derivatives satisfy
  \[
    \|\partial_t\widehat m_p\|_{L^2_tL^2_x}
      +\|\nabla\widehat m_p\|_{L^2_tL^2_x}
    \le\gamma_p^{-1}\bigl(
       \|\partial_t\widehat q^{\,*}\|_{L^2_tL^2_x}
       +\|\nabla\widehat q^{\,*}\|_{L^2_tL^2_x}\bigr).
  \]
  Thus no separate control or active-set regularity assumption is
  needed. 
\end{rmk}

\begin{proof}
  Write $m=\widehat m_p$, $m_h=\widehat m_{p,h}$. 
  We use the
  $L^2((0,T)\times\Omega)$ inner product; discrete state and adjoint
  quantities are identified with their piecewise-constant-in-time
  extensions.  The discrete functional is evaluated at controls in
  $\mathcal M_{p,ad}^h$.

  By the discrete strong monotonicity
  (Lemma~\ref{lem:tp:Jh-monotone}) and the discrete variational
  inequality at $m_h$ tested with $\Pi_h^M m$, setting
  $e_m^h=m_h-\Pi_h^M m$ gives
  \begin{align}
    \label{eq:tp:main-start}
    \gamma_p\|e_m^h\|_{L^2_tL^2_x}^2
    &\le (j_h'(m_h)-j_h'(\Pi_h^M m),e_m^h)_{L^2_tL^2_x}
     \le - (j_h'(\Pi_h^M m),e_m^h)_{L^2_tL^2_x}.
  \end{align}
  Add and subtract the continuous reduced gradient
  $j'(m)=\gamma_pm+\widehat q^{\,*}$ (Lemma~\ref{lem:tp:VI}).  Because $m_h$ is admissible,
  the continuous variational inequality gives
  $\langle j'(m),m_h-m\rangle\ge0$.  Hence
  \begin{align}
    \label{eq:tp:after-continuous-VI}
    \gamma_p\|e_m^h\|_{L^2_tL^2_x}^2
    &\le
      -(j_h'(\Pi_h^M m)-j'(m),e_m^h)_{L^2_tL^2_x}
      +(j'(m),\Pi_h^M m-m)_{L^2_tL^2_x}\notag\\
    &\le
      \bigl| (j_h'(\Pi_h^M m)-j'(m),e_m^h)_{L^2_tL^2_x} \bigr|
      +(j'(m),\Pi_h^M m-m)_{L^2_tL^2_x}.
  \end{align}
  The $L^2$ orthogonality of the space-time projection and
  $\Pi_h^M m=\Pi_h^M m$ give
  \begin{align*}
    (j'(m),\Pi_h^M m-m)_{L^2_tL^2_x}
      &=-((I-\Pi_h^M)j'(m),(I-\Pi_h^M)m)_{L^2_tL^2_x}
       \le C(h+\Delta t)^2,
  \end{align*}
  by Assumption~R and the derived regularity
  \eqref{eq:tp:control-regularity}.

  Since $e_m^h\in M_{p,h}$, only the $M_{p,h}$-projection of
  the discrete--continuous gradient difference matters.  Using
  \eqref{eq:tp:Jh-gradient} and the continuous gradient identity,
  and observing that $\Pi_h^M(\Pi_h^M m-m)=0$, we obtain
  \begin{align*}
    \Pi_h^M\bigl(j_h'(\Pi_h^M m)-j'(m)\bigr) =\Pi_h^M\bigl(q_h^*(\Pi_h^M m)-q_h^*(m)\bigr)
          +\Pi_h^M\bigl(q_h^*(m)-\widehat q^{\,*}\bigr).
  \end{align*}
  The control-to-adjoint sensitivity estimate
  \eqref{eq:tp:diff-adj-est}, Jensen's inequality, and the
  projection error of $m$ imply
  \[
    \|q_h^*(\Pi_h^M m)-q_h^*(m)\|_{L^2_tL^2_x}
        \le C\|\Pi_h^M m-m\|_{L^2_tL^2_x}
        \le C(h+\Delta t).
  \]
  The fixed-control adjoint error
  \eqref{eq:tp:adjoint-error-fixed-control}, together with
  the time-averaging projection error of $\widehat q^{\,*}$, gives
  \[
    \|\Pi_h^M(q_h^*(m)-\widehat q^{\,*})\|_{L^2_tL^2_x}
       \le C(h+\Delta t).
  \]
  Therefore
  \[
    \bigl|( j_h'(\Pi_h^M m)-j'(m),e_m^h)_{L^2_tL^2_x}\bigr|
      \le C(h+\Delta t)\|e_m^h\|_{L^2_tL^2_x}.
  \]
  Substitution into \eqref{eq:tp:after-continuous-VI} gives
  \[
    \gamma_p\|e_m^h\|^2
       \le C(h+\Delta t)\|e_m^h\|_{L^2_tL^2_x}+C(h+\Delta t)^2,
  \]
  and hence $\|e_m^h\|\le C(h+\Delta t)$.
  Finally,
  \[
    \|m-m_h\|_{L^2_tL^2_x}\le\|m-\Pi_h^Mm\|_{L^2_tL^2_x}+\|e_m^h\|_{L^2_tL^2_x}
      \le C(h+\Delta t),
  \]
  proving \eqref{eq:tp:control-error}.
\end{proof}

\begin{theorem}
  \label{thm:tp:state-error}
  Under the assumptions of Theorem~\ref{thm:tp:control-error}
  and the regularity \eqref{eq:tp:reg-state-adj},
  the state error satisfies
  \begin{align}
    \label{eq:tp:state-error}
    \Bigl(\Delta t\sum_{j=1}^N
    \bigl(\|\overline{\widehat u}{}^{\,j} - \widehat{u}_h^j\|_1^2
    + \|\overline{\widehat p_T}{}^{\,j} - \widehat{p}_{T,h}^j\|_0^2
    + \|\overline{\widehat p}{}^{\,j} - \widehat{p}_h^j\|_0^2\bigr)
    \Bigr)^{1/2}
    \le C(h+\Delta t),
  \end{align}
  where $C>0$ is independent of $h$ and $\Delta t$, uniform for $\lambda\ge\lambda_0>0$.
\end{theorem}

\begin{proof}
  We split
  \begin{align*}
    \overline{\widehat u}{}^{\,j} - \widehat{u}_h^j
    = \bigl(\overline{\widehat u}{}^{\,j} - u_h^j(\widehat{m}_p)\bigr)
    + \bigl(u_h^j(\widehat{m}_p) - \widehat{u}_h^j\bigr),
  \end{align*}
  and analogously for $\overline{\widehat p_T}{}^{\,j}$ and $\overline{\widehat p}{}^{\,j}$.
  The first term is the fixed-control discretization error, bounded
  by $O(h+\Delta t)$ in the $\ell^2_t$ sense by both components of
  Lemma~\ref{lem:tp:fixed-control-state-error}: \eqref{eq:tp:fixed-state-mech}
  for the $(u,p_T)$-part and \eqref{eq:tp:fixed-state-diff} for
  the $p$-part.
  Together with Lemma~\ref{lem:tp:generic-stability}, one obtains
  \begin{align}
    \label{eq:tp:fixed-control-state-error}
    \Bigl(\Delta t\sum_{j=1}^N\bigl(
      \|\overline{\widehat u}{}^{\,j} - u_h^j(\widehat m_p)\|_1^2
      + \|\overline{\widehat p_T}{}^{\,j} - p_{T,h}^j(\widehat m_p)\|_0^2
      + \|\overline{\widehat p}{}^{\,j} - p_h^j(\widehat m_p)\|_0^2
    \bigr)\Bigr)^{1/2}
    &\le C(h+\Delta t).
  \end{align}
  The second term is the control-perturbation error; by the
  $\ell^2_t$ estimates of Lemma~\ref{lem:tp:control-sensitivity},
  \begin{align*}
    \Bigl(\Delta t\sum_{j=1}^N\bigl(
    \|u_h^j(\widehat m_p) - \widehat{u}_h^j\|_V^2
    + \|p_{T,h}^j(\widehat m_p) - \widehat{p}_{T,h}^j\|_0^2
    + \|p_h^j(\widehat m_p) - \widehat{p}_h^j\|_0^2
    \bigr)\Bigr)^{1/2}
    \le C(h+\Delta t),
  \end{align*}
  where the last step uses Theorem~\ref{thm:tp:control-error}.
  Combining with \eqref{eq:tp:fixed-control-state-error} gives
  \eqref{eq:tp:state-error}.
\end{proof}

\begin{rmk}
  \label{rmk:tp:rate}
  The fixed-control state estimate of
  Lemma~\ref{lem:tp:fixed-control-state-error} achieves the rate
  $O(h+\Delta t)$.  This follows because, with time-independent
  coefficients and the cell-average formulation, the stiffness residual
  $R_S^{k+1}(q_h)=(\kappa_p\nabla(I-P_h)\bar p^{k+1},\nabla q_h)_\Omega$ is purely spatial, giving
  $O(h)$.  The mass residual and the commutator are estimated jointly:
  after summation by parts they contribute $O(h)$ through
  Lemma~\ref{lem:tp:mass-consistency} and $O(\Delta t)$ through the
  cell-average defect $d^n$.

  The final control and optimal-state estimates
  (Theorems~\ref{thm:tp:control-error}--\ref{thm:tp:state-error}) achieve
  the same rate $O(h+\Delta t)$.  This follows because the adjoint
  fixed-control estimate of
  Lemma~\ref{lem:tp:adjoint-error-fixed-control} is controlled entirely
  through the $O(h+\Delta t)$ fixed-control state estimate,
  applied both to the forward-in-time system and, after time reversal,
  to the reduced scalar adjoint system, which has the identical symmetric structure.
\end{rmk}

\section{Numerical experiments}
\label{sec:numerical}
\newcommand{\numObjectiveReduction}{98.72}

We examine convergence of the fully discrete optimality system,
accuracy as the material parameters vary, and a localized pressure-control
problem. All results in this section were computed with NGSolve 6.2.2608
\cite{schoeberl2014ngsolve}. The error estimates are uniform in $\lambda$
for storage bounded below by a positive constant, while the
Schur-complement formulation remains well posed at $s_0=0$ for each
fixed finite $\lambda$. The experiments examine both regimes,
including parameter combinations for which the analysis does not
establish simultaneous uniformity. The implementation, input
configurations, and raw results are supplied in the accompanying
supplementary code archive.

Throughout, $\Omega=(0,1)^2$, $T=1$, $\mu=\alpha=1$, and
$\kappa_p=0.1$. We impose $u=0$ on $\Gamma_D=\{0\}\times(0,1)$,
zero traction on the remaining three sides, and $p=0$ on the whole boundary.
Thus $|\Gamma_N|>0$ and the total-pressure space is the full continuous
$\mathbb P_1$ space, without a mean constraint or pressure pinning.
The dG(0) time treatment of bounded controls is related to the
fully discrete optimal-control analysis in \cite{MR3022219}.
Each square of a uniform $n\times n$ grid is divided into two triangles,
so that $h=\sqrt{2}/n$. The state and adjoint use the
$\mathbb P_2/\mathbb P_1/\mathbb P_1$ spaces of
Section~\ref{sec:THM:FE-spaces}; controls are constant on each triangle
and time slab. We vary $n$ and $N$ independently in the convergence tests;
each table states both resolutions. The localized example uses
$n=N=16$ and $32$.
Both tracking weights are one.

\subsection{Convergence verification}
\label{sec:num:convergence}

We use a manufactured optimum that requires no additional fluid source.
Define
\begin{align*}
 B(x,y)&=4096x^3(1-x)^3y^3(1-y)^3,
 &w&=\tfrac1{20}(\partial_yB,-\partial_xB),\\
 \eta(t)&=t^2(1-t)^3,&\beta(t)&=(1-t)^2.
\end{align*}
With $\gamma_p=0.1$, the exact state, control, and adjoint are
\begin{align}
 u&=\eta w,&p&=\eta B,&p_T&=\alpha p,\notag\\
 m_p&=s_0\eta' B-\kappa_p\eta\Delta B,
 &q^*&=-\gamma_pm_p,&q_T^*&=-\alpha q^*,\notag\\
 v^*&=\beta w.&&&
 \label{eq:num:mms}
\end{align}
Since $\operatorname{div}w=0$, the total-pressure relations hold for
all $\lambda$. The bubble and its derivatives ensure the homogeneous
boundary conditions, including zero traction on $\Gamma_N$.
The initial state and terminal adjoint are zero. The mechanical forcing
in the differentiated state equation and the targets are set to
\begin{gather*}
 f=-\operatorname{div}(2\mu\varepsilon(\dot u))+\nabla\dot p_T, \qquad \qquad 
  R_u=\operatorname{div}(2\mu\varepsilon(\dot v^*))-\nabla\dot q_T^*,
 \\
 u_C=u-R_u, \qquad  R_p=-s_0\dot q^*-\kappa_p\Delta q^*, \qquad p_C=p-R_p.
\end{gather*}
The two $\alpha^2/\lambda$ terms in the scalar adjoint equation cancel
for this family. Thus the state equations, adjoint equations, and
$\gamma_pm_p+q^*=0$ hold simultaneously. Strict convexity identifies
\eqref{eq:num:mms} as the unique optimum. We solve the unconstrained
problem here, equivalently an inactive sufficiently wide box problem.
The exact fields and targets remain bounded independently of $\lambda$.
This divergence-free family isolates parameter effects but does not
sample general volumetric deformation.

For state variables, the reported errors are
\begin{align*}
 E_u^2&=\Delta t\sum_k\|\bar{u}^k -u_h^k\|_{H^1}^2,
 &E_T^2&=\Delta t\sum_k\|\bar{p_T}^k-p_{T,h}^k\|_0^2,\\
 E_p^2&=\Delta t\sum_k\|\bar{p}^k -p_h^k\|_0^2,
 &E_m&=\|m_p-m_{p,h}\|_{L^2_tL^2_x}.
\end{align*}
We set $E_{\rm state}=(E_u^2+E_T^2+E_p^2)^{1/2}$.
The adjoint errors $E_{v^*}$, $E_{q_T^*}$, and $E_{q^*}$ are true
space--time norms in $H^1$, $L^2$, and $L^2$, respectively, with
$z_h^k$ held constant on $(t_k,t_{k+1}]$. In particular, both the control
and adjoint errors include the within-slab variation of the exact solution.
The coupled refinements $n=N=8,16,32,64$ are retained in the
supplementary data. Here we separate the spatial and temporal refinements
to examine the contributions to the bounds in
Theorems~\ref{thm:tp:control-error} and~\ref{thm:tp:state-error}.

\emph{Spatial convergence.}
To separate the spatial and temporal effects, we keep $N=2048$
fixed and use $n=8,16,32$, with $\lambda=s_0=1$. All errors in
Table~\ref{tab:num:spatial-independent} are measured against the analytic
optimum, with the same state cell-average and true space--time norms as
above. Spatial rates are
$r_x=\log_2(E(n,N)/E(2n,N))$ with $N$ unchanged.
The time-resolution check repeats the finest spatial level with $N=4096$.
The measured changes in all seven errors are reported below together
with the parameter study; they check residual temporal contamination
rather than assume it is negligible from the size of $\Delta t$ alone.

\begin{table}[htbp]
\centering\scriptsize
\caption{Spatial refinement at fixed $N=2048$, $\lambda=s_0=1$: errors relative to the manufactured optimum. Rates in parentheses compare consecutive spatial levels.}\label{tab:num:spatial-independent}
\begin{tabular}{rrrrrrrrrrr}\hline
$n$ & $N$ & $h$ & $\Delta t$ & $E_u$ & $E_T$ & $E_p$ & $E_m$ & $E_{v^*}$ & $E_{q_T^*}$ & $E_{q^*}$\\\hline
8 & 2048 & 1.768e-01 & 4.883e-04 & \shortstack{1.256e-03\\(--)} & \shortstack{6.197e-04\\(--)} & \shortstack{9.201e-04\\(--)} & \shortstack{1.106e-02\\(--)} & \shortstack{2.601e-02\\(--)} & \shortstack{5.592e-04\\(--)} & \shortstack{6.140e-04\\(--)}\\[3pt]
16 & 2048 & 8.839e-02 & 4.883e-04 & \shortstack{3.342e-04\\(1.91)} & \shortstack{1.756e-04\\(1.82)} & \shortstack{2.636e-04\\(1.80)} & \shortstack{5.199e-03\\(1.09)} & \shortstack{6.874e-03\\(1.92)} & \shortstack{1.193e-04\\(2.23)} & \shortstack{1.614e-04\\(1.93)}\\[3pt]
32 & 2048 & 4.419e-02 & 4.883e-04 & \shortstack{8.615e-05\\(1.96)} & \shortstack{4.630e-05\\(1.92)} & \shortstack{6.929e-05\\(1.93)} & \shortstack{2.547e-03\\(1.03)} & \shortstack{1.759e-03\\(1.97)} & \shortstack{2.976e-05\\(2.00)} & \shortstack{4.179e-05\\(1.95)}\\[3pt]
\hline\end{tabular}
\end{table}

\emph{Temporal convergence.}
At fixed $n=32$, we use $N=8,16,32,64,128$. Spatial discretization
produces an error floor in the comparison with the analytic optimum.
For example, a separate $n=64$ pilot gives analytic control errors
$2.270\times10^{-3}$ and $1.582\times10^{-3}$ at $N=64$ and $128$,
respectively, an apparent temporal rate of only $0.52$.
To measure temporal error separately, we compare each solution with a
reference optimum on the \emph{same} spatial mesh, computed with
$N_{\rm ref}=4096$. These reference errors approximate errors relative
to the continuous-in-time, spatially discrete optimum; they are distinct
from the total errors against \eqref{eq:num:mms}.
For $r=N_{\rm ref}/N$, the displacement reference error is
\[
 (E_u^t)^2=\Delta t\sum_{k=1}^{N}
 \left\|u_{h,N}^{k}
 -\frac1r\sum_{j=(k-1)r+1}^{kr}u_{h,N_{\rm ref}}^{j}\right\|_{H^1}^2.
\]
The other state reference errors use the corresponding $L^2$ norms.
Control and adjoint reference errors integrate the difference of their
piecewise-constant reconstructions on the common fine time grid, retaining
variation within the coarse slabs. For example,
$E_m^t=\|m_{h,N}-m_{h,N_{\rm ref}}\|_{L^2_tL^2_x}$.
Rates are $r_t=\log_2(E^t(N)/E^t(2N))$ at fixed $n$ and $N_{\rm ref}$.

\begin{table}[htbp]
\centering\scriptsize
\caption{Temporal refinement at fixed $n=32$, $\lambda=s_0=1$: reference errors $E^t$ using $N_{\rm ref}=4096$ on the same spatial mesh. Rates are in parentheses. These quantities approximate the temporal discretization errors relative to the continuous-in-time, spatially discrete optimum.}\label{tab:num:temporal-independent}
\begin{tabular}{rrrrrrrrrrr}\hline
$n$ & $N$ & $h$ & $\Delta t$ & $E_u^t$ & $E_T^t$ & $E_p^t$ & $E_m^t$ & $E_{v^*}^t$ & $E_{q_T^*}^t$ & $E_{q^*}^t$\\\hline
32 & 8 & 4.419e-02 & 1.250e-01 & \shortstack{3.679e-03\\(--)} & \shortstack{9.053e-04\\(--)} & \shortstack{7.279e-04\\(--)} & \shortstack{1.443e-02\\(--)} & \shortstack{6.303e-02\\(--)} & \shortstack{1.516e-03\\(--)} & \shortstack{1.445e-03\\(--)}\\[3pt]
32 & 16 & 4.419e-02 & 6.250e-02 & \shortstack{1.855e-03\\(0.99)} & \shortstack{4.713e-04\\(0.94)} & \shortstack{3.909e-04\\(0.90)} & \shortstack{7.391e-03\\(0.97)} & \shortstack{3.114e-02\\(1.02)} & \shortstack{7.643e-04\\(0.99)} & \shortstack{7.398e-04\\(0.97)}\\[3pt]
32 & 32 & 4.419e-02 & 3.125e-02 & \shortstack{9.269e-04\\(1.00)} & \shortstack{2.397e-04\\(0.98)} & \shortstack{2.031e-04\\(0.94)} & \shortstack{3.707e-03\\(1.00)} & \shortstack{1.544e-02\\(1.01)} & \shortstack{3.804e-04\\(1.01)} & \shortstack{3.711e-04\\(1.00)}\\[3pt]
32 & 64 & 4.419e-02 & 1.562e-02 & \shortstack{4.603e-04\\(1.01)} & \shortstack{1.201e-04\\(1.00)} & \shortstack{1.031e-04\\(0.98)} & \shortstack{1.846e-03\\(1.01)} & \shortstack{7.653e-03\\(1.01)} & \shortstack{1.887e-04\\(1.01)} & \shortstack{1.848e-04\\(1.01)}\\[3pt]
32 & 128 & 4.419e-02 & 7.812e-03 & \shortstack{2.265e-04\\(1.02)} & \shortstack{5.940e-05\\(1.02)} & \shortstack{5.133e-05\\(1.01)} & \shortstack{9.131e-04\\(1.02)} & \shortstack{3.776e-03\\(1.02)} & \shortstack{9.314e-05\\(1.02)} & \shortstack{9.140e-05\\(1.02)}\\[3pt]
\hline\end{tabular}
\end{table}

\begin{figure}[htbp]
 \centering
 \includegraphics[width=\textwidth]{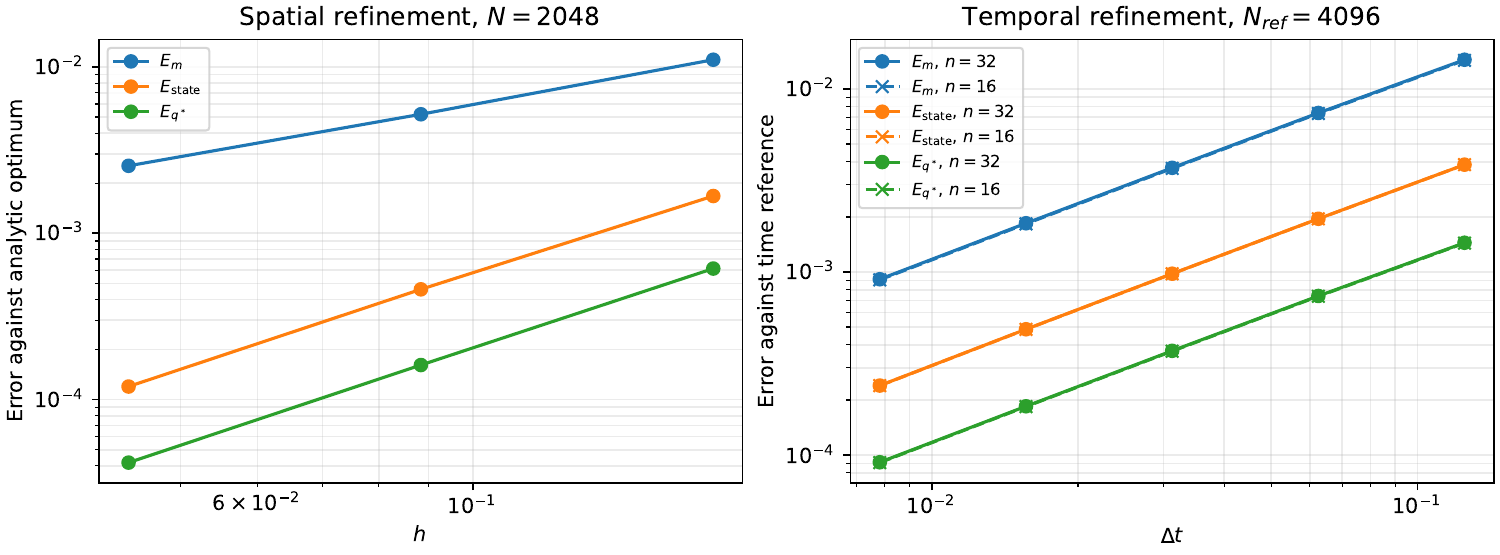}
 \caption{Independent refinements at $\lambda=s_0=1$.
 Left: errors against the analytic optimum at fixed $N=2048$.
 Right: temporal reference errors at fixed $n=16$ and $32$, each using
 its own same-mesh $N_{\rm ref}=4096$ optimum.}
 \label{fig:num:independent}
\end{figure}

We repeat the reference comparison with $N_{\rm ref}=2048$ to quantify
reference sensitivity. We also repeat the temporal refinement on the
$n=16$ mesh, with references computed on that same mesh. Figure~\ref{fig:num:independent}
compares both temporal sequences, so that the observed temporal orders
can be assessed on two fixed spatial discretizations. The analytic-error
and reference-error data are both included in the supplementary results.

\subsection{Robustness with respect to \texorpdfstring{$\lambda$ and $s_0$}{lambda and s0}}
\label{sec:num:robustness}

We next solve the manufactured optimal control problem for all pairs
\[
 \lambda\in\{1,10^2,10^4,10^6\},\qquad
 s_0\in\{1,10^{-2},10^{-4},0\},
\]
using the same spatial sequence $n=8,16,32$ at fixed $N=2048$
and temporal sequence $N=8,16,32,64,128$ at fixed $n=32$.
Each parameter pair has its own same-mesh $N_{\rm ref}=4096$ reference
optimum, with the $N_{\rm ref}=2048$ comparison described above.
The data are reconstructed from \eqref{eq:num:mms} for each $s_0$;
for fixed $s_0$, neither the exact solution nor the targets change
with $\lambda$.

For each fixed positive storage coefficient, comparison across
$\lambda$ tests the nearly incompressible regime covered by the
$\lambda$-uniform error analysis. Comparison at fixed $\lambda$
tests the approach to zero storage. The runs with $s_0=0$ use the
same matrices, boundary conditions, and optimization procedure:
no positive storage floor or pressure penalty is inserted.
Table~\ref{tab:num:independent-robustness} reports absolute error
magnitudes and the final spatial and temporal rates separately.
Here the superscript $x$ denotes the analytic error at fixed fine
$\Delta t$, while the superscript $t$ denotes the temporal reference
error at fixed $h$. The combined state error uses the three state
norms defined above.

\begin{table}[htbp]
\centering\scriptsize
\caption{Independent parameter study. Spatial errors use $(n,N)=(32,2048)$; spatial rates compare $n=16,32$ at fixed $N=2048$. Temporal reference errors use $(n,N)=(32,128)$ and $N_{\rm ref}=4096$; temporal rates compare $N=64,128$ at fixed $n=32$.}\label{tab:num:independent-robustness}
\begin{tabular}{rrrrrrrrrr}\hline
$\lambda$ & $s_0$ & $E_m^x$ & $r_m^x$ & $E_{\rm state}^x$ & $r_{\rm state}^x$ & $E_m^t$ & $r_m^t$ & $E_{\rm state}^t$ & $r_{\rm state}^t$\\\hline
1 & 1 & 2.547e-03 & 1.03 & 1.199e-04 & 1.94 & 9.131e-04 & 1.02 & 2.398e-04 & 1.02\\
1 & 0.01 & 2.287e-03 & 1.04 & 1.277e-04 & 1.94 & 2.835e-04 & 1.00 & 2.408e-04 & 1.02\\
1 & 0.0001 & 2.287e-03 & 1.04 & 1.278e-04 & 1.94 & 2.820e-04 & 1.00 & 2.407e-04 & 1.02\\
1 & 0 & 2.287e-03 & 1.04 & 1.278e-04 & 1.94 & 2.820e-04 & 1.00 & 2.407e-04 & 1.02\\
100 & 1 & 2.548e-03 & 1.03 & 1.108e-04 & 1.94 & 8.741e-04 & 1.01 & 2.507e-04 & 1.02\\
100 & 0.01 & 2.287e-03 & 1.04 & 1.142e-04 & 1.95 & 2.615e-04 & 1.00 & 2.454e-04 & 1.02\\
100 & 0.0001 & 2.287e-03 & 1.04 & 1.142e-04 & 1.95 & 2.614e-04 & 1.00 & 2.454e-04 & 1.02\\
100 & 0 & 2.287e-03 & 1.04 & 1.142e-04 & 1.95 & 2.614e-04 & 1.00 & 2.454e-04 & 1.02\\
10000 & 1 & 2.548e-03 & 1.03 & 1.108e-04 & 1.94 & 8.725e-04 & 1.01 & 2.514e-04 & 1.02\\
10000 & 0.01 & 2.287e-03 & 1.04 & 1.142e-04 & 1.95 & 2.615e-04 & 1.00 & 2.460e-04 & 1.02\\
10000 & 0.0001 & 2.287e-03 & 1.04 & 1.142e-04 & 1.95 & 2.614e-04 & 1.00 & 2.460e-04 & 1.02\\
10000 & 0 & 2.287e-03 & 1.04 & 1.142e-04 & 1.95 & 2.614e-04 & 1.00 & 2.460e-04 & 1.02\\
1e+06 & 1 & 2.548e-03 & 1.03 & 1.108e-04 & 1.94 & 8.724e-04 & 1.01 & 2.514e-04 & 1.02\\
1e+06 & 0.01 & 2.287e-03 & 1.04 & 1.142e-04 & 1.95 & 2.615e-04 & 1.00 & 2.460e-04 & 1.02\\
1e+06 & 0.0001 & 2.287e-03 & 1.04 & 1.142e-04 & 1.95 & 2.614e-04 & 1.00 & 2.460e-04 & 1.02\\
1e+06 & 0 & 2.287e-03 & 1.04 & 1.142e-04 & 1.95 & 2.614e-04 & 1.00 & 2.460e-04 & 1.02\\
\hline\end{tabular}
\end{table}

The final spatial control rates across the parameter grid are $1.03\text{--}1.04$, and the combined state rates are $1.94\text{--}1.95$. The smooth state displays higher spatial orders than the first-order control, consistent with the different approximation spaces; this does not strengthen the theorem for minimally regular solutions. The final temporal reference rates are $1.00\text{--}1.02$ for the control and $1.02$ for the combined state.

For $\lambda=s_0=1$, replacing $N=2048$ by $4096$ at $n=32$ changes the seven analytic errors by at most 2.72\%, while doubling the temporal reference changes the temporal errors by at most 3.26\%. Across the parameter grid, at fixed $n=32$, replacing $N=2048$ by $4096$ changes the analytic errors by at most 10.17\% across all seven variables and sixteen parameter pairs. Replacing the temporal reference $N_{\rm ref}=2048$ by $4096$ changes the reported temporal errors by at most 3.26\% over all five temporal levels and all parameter pairs. These are measured resolution-sensitivity checks, not rigorous bounds on the remaining reference error.

For the control and combined state quantities displayed in Table~\ref{tab:num:independent-robustness}, the respective spatial and temporal sensitivity maxima are 1.20\% and 3.23\%. The largest individual-variable sensitivity concerns total pressure in the nearly incompressible cases. All component sensitivities are retained in the supplementary data; we do not infer equally resolved individual-variable rates for every parameter pair.

The nearly incompressible and zero-storage cases retain the observed control and combined-state orders in this manufactured family.

A separate algebraic check makes the limitation of this observation
explicit. On the $n=8$ mesh, we assemble the discrete effective-mass matrix corresponding to $m_{\mathrm{eff},h}$ in \eqref{eq:tp:discrete-mass-def} by eliminating the displacement and the total pressure variables, and compute its smallest generalized eigenvalue relative to the fluid-pressure $L^2$ mass matrix. At $s_0=0$, its smallest
eigenvalue is approximately $2.9595\times10^{-1}$,
$9.7625\times10^{-3}$, $9.9976\times10^{-5}$, and
$1.0000\times10^{-6}$ for the four values of $\lambda$ above.
The effective mass is positive for every tested finite $\lambda$, while
its coercivity decreases as $\lambda$ grows. Consequently, the bounded
errors in this smooth, divergence-free manufactured family do not imply
uniform coercivity, or a uniform error theorem, in the simultaneous
limit $s_0\downarrow0$, $\lambda\to\infty$.

\subsection{A localized pressure-control example}
\label{sec:num:localized}

Finally, we consider a problem without a prescribed analytic optimum.
The body force and initial state vanish, $u_C=0$, and
\[
 p_C(t,x,y)=0.15\,t(2-t)\sin(\pi x)\sin(\pi y).
\]
The source is restricted to the mesh-aligned actuator
$\omega=(1/4,3/4)^2$: admissible controls satisfy $0\le m_p\le0.4$
on $(0,1)\times\omega$ and $m_p=0$ outside it.
Equivalently, the distributed source is extended by zero from $\omega$;
the same discrete projection uses equal zero bounds on the other cells.
This spatially restricted admissible set extends the constant-box
setting of Section~\ref{sec:OCP}. The control may jump across
$\partial\omega$, so this example illustrates the discrete optimality
system without assuming global $H^1(\Omega)$ control regularity or
claiming the rate of Theorem~\ref{thm:tp:control-error}.
We choose $\lambda=10^4$, $s_0=0$, and $\gamma_p=10^{-3}$.
All remaining coefficients and boundary conditions are unchanged.
This finite-parameter example combines nearly incompressible mechanics,
zero storage, and an active upper control bound.

\begin{table}[htbp]
\centering\small
\caption{Localized source control. The tracking norm includes displacement and pressure. Active fractions refer to the whole space--time domain and exclude fixed-zero cells.}\label{tab:num:localized}
\begin{tabular}{rrrrrr}\hline
$n=N$ & $J_h$ (continuous target) & tracking norm & $\|m_h\|$ & upper active & $r_{\rm proj}$ \\\hline
16 & 1.985e-05 & 4.425e-03 & 1.418e-01 & 7.86\% & 9.527e-13 \\
32 & 1.926e-05 & 4.318e-03 & 1.410e-01 & 7.76\% & 6.887e-13 \\
\hline\end{tabular}
\end{table}

Table~\ref{tab:num:localized} reports the objective, tracking norm,
control norm, and upper-active fraction for two mesh/time resolutions.
The reported objective uses the continuous target norm; it differs from
the slab-average-target objective in \eqref{eq:tp:Jh} only by a known,
control-independent constant. The two objectives therefore have the
same minimizer and gradient. With zero control, the state is zero,
which supplies the uncontrolled baseline $J(0)=1.5\times10^{-3}$.
At $n=N=32$, the objective is reduced by \numObjectiveReduction\%.
The optimized source reduces the tracking discrepancy, while its limited
support, the upper bound, and regularization prevent exact matching.
The active fraction is measured relative to the whole space--time domain;
fixed-zero cells outside the actuator are excluded from the active-set
statistics. 

\begin{figure}[htbp]
 \centering
 \includegraphics[width=\textwidth]{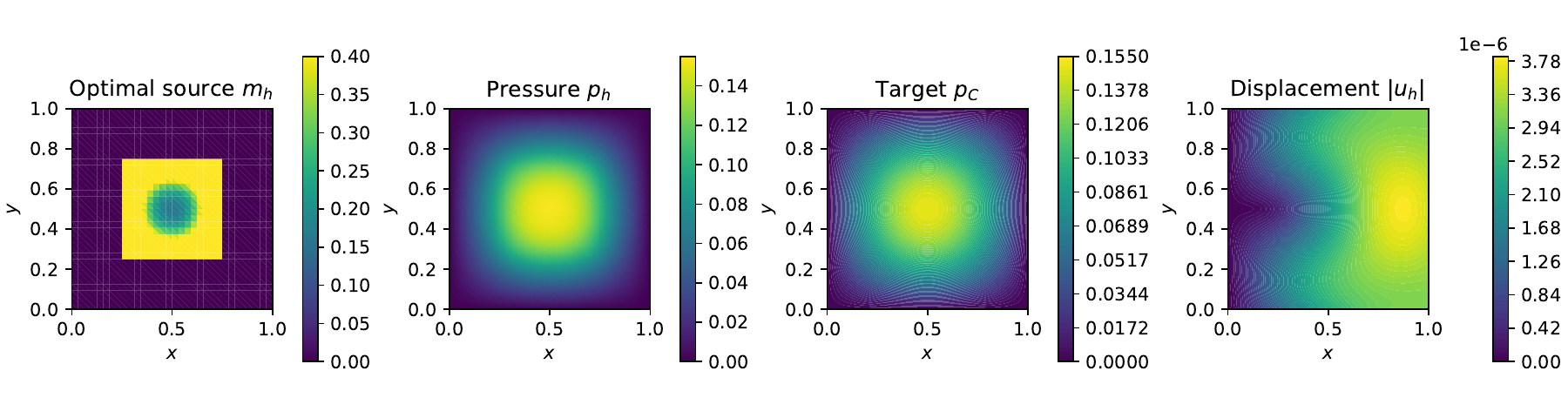}
 \caption{Localized control at $n=N=32$: last-slab source, terminal
 pressure, target at $T=1$, and displacement magnitude. The source is
 shown with its exact constant value on each triangle, so individual
 elements remain visible near the center. The displacement is resolved
 with its own color scale because it is small in this nearly
 incompressible regime.}
 \label{fig:num:localized}
\end{figure}


\appendix

\section{Regularity of the optimal state and adjoint}
\label{app:regularity}

This appendix gives sufficient conditions for Assumption~R.

The regularity package established below is
\begin{align}
  \label{eq:app:reg-state}
  \widehat u,\,\widehat v^*
  &\in L^2(0,T;\boldsymbol{H}^2(\Omega))\cap H^1(0,T;V),
  \\
  \label{eq:app:reg-pT}
  \widehat p_T,\,\widehat q_T^*
  &\in L^2(0,T;H^1(\Omega))\cap H^1(0,T;Q_T),
  \\
  \label{eq:app:reg-diff}
  \widehat p,\,\widehat q^{\,*}
  &\in L^2(0,T;H^2(\Omega))\cap H^1(0,T;L^2(\Omega)).
\end{align}
Here $Q_T=L^2(\Omega)$.  In particular,
\eqref{eq:app:reg-diff} does \emph{not} assert
$\widehat p,\widehat q^{\,*}\in H^1(0,T;Q)$;
such a statement would require additional compatibility assumptions
on the time-differentiated system and is not needed for the dG(0)
error estimate.

The basic state and adjoint time-regularity below follows from the
Schur complement estimate for every $s_0\ge0$.  When these regularity bounds are combined with the
$\lambda$-uniform error estimates of Section~\ref{sec:THM:error}, we
retain the positive-storage regime stated there.

We also assume the following scalar elliptic and Stokes-type regularity
estimates.
\begin{enumerate}
\item \label{item:app:elliptic} \emph{Scalar $H^2$-regularity.}  For every $F\in L^2(\Omega)$,
  the solution $w\in H^1_0(\Omega)$ of
  $-\operatorname{div}(\kappa\nabla w)=F$ satisfies
  $\|w\|_{H^2(\Omega)}\le C_{\rm ell}\|F\|_{L^2(\Omega)}$.
  This holds for convex polygonal/polyhedral $\Omega$ or for
  $\partial\Omega\in C^{1,1}$, for constant $\kappa>0$ or for
  sufficiently smooth, uniformly elliptic coefficients.  It is applied
  with $\kappa=\kappa_p$.
\item \label{item:app:stokes} \emph{Stokes-type $H^2\times H^1$-regularity.}  If
  $(w,r)\in V\times Q_T$ solves
  $a_u(w,v)+b(r,v)=(F,v)_\Omega$ and
  $b(q_T,w)-\lambda^{-1}(r,q_T)_\Omega=(\zeta,q_T)_\Omega$
  for all $(v,q_T)\in V\times Q_T$, then
  $\|w\|_{\boldsymbol{H}^2(\Omega)}+\|r\|_{H^1(\Omega)}
  \le C_{\rm St}(\|F\|_{L^2(\Omega)}+\|\zeta\|_{H^1(\Omega)})$.
  We assume this estimate holds for the particular domain $\Omega$,
  boundary conditions, and Lam\'e parameters under consideration, with
  $C_{\rm St}$ independent of $\lambda$ for $\lambda\ge\lambda_0>0$;
  whether it holds on a given convex polyhedral domain with mixed
  boundary conditions depends on additional geometric and compatibility
  conditions.  Note that $\zeta$ enters through an $H^1$-norm because
  it plays the role of a divergence datum; $L^2$-control of $\zeta$
  alone is not sufficient for $H^1$-regularity of $r$.
\end{enumerate}

The argument proceeds in three steps:
(i)~basic energy regularity of the state, giving in particular
$\dot p\in L^2(0,T;L^2(\Omega))$ and $p\in L^\infty(0,T;Q)$
(Proposition~\ref{prop:app:basic-state-reg});
(ii)~spatial $H^2\times H^1$-regularity of $(u,p_T)$ via the
Stokes-type estimate applied to the integrated mechanical equation
(Proposition~\ref{prop:app:stokes-reg}); and
(iii)~spatial $H^2$-regularity of $p$ via the scalar
elliptic estimate applied to the diffusion equation
(Proposition~\ref{prop:app:elliptic-reg}).
The adjoint is handled in Proposition~\ref{prop:app:adjoint-reg}
by time-reversal, applying steps (i)--(iii) to the reversed system
and then reading off the spatial regularity.

\subsection{Basic energy regularity of the state}
\label{app:basic-state-reg}

\begin{proposition}
  \label{prop:app:basic-state-reg}
  Let $f:=\dot{\tilde f}$.  Assume
  \begin{align}
    \label{eq:app:data-energy}
    f \in L^2(0,T;L^2(\Omega)^d),
    \qquad
    g_p \in L^2(0,T;L^2(\Omega)),
  \end{align}
  and compatible initial data $(u_0,p_0)\in V\times Q$
  for the total pressure formulation.
  Let $(u,p_T,p)$ be the weak solution of
  \eqref{eq:tp:strong-pT}.  Then
  \begin{align}
    \label{eq:app:basic-reg-est}
    u &\in H^1(0,T;V),
    \qquad
    p_T \in H^1(0,T;Q_T),
    \notag\\
    p &\in L^\infty(0,T;Q)\cap H^1(0,T;L^2(\Omega)),
  \end{align}
  with
  \begin{align}
    \label{eq:app:basic-reg-bound}
    \|\dot u\|_{L^2_tV}
    + \|\dot p_T\|_{L^2_tQ_T}
    + \|\dot p\|_{L^2_tL^2_x}
    + \|p\|_{L^\infty_tQ}
    \le C_{\rm en},
  \end{align}
  where $C_{\rm en}$ depends on the data and material parameters,
  is independent of $h$, $\Delta t$, and may depend on $c_{\rm eff}$ in \eqref{eq:tp:reduced-mass-identity}. Dependence on
  $\kappa_p$ is through the weighted norms and elliptic regularity constants.
\end{proposition}

\begin{proof}
  This is exactly the storage-robust estimate of
  Lemma~\ref{lem:tp:energy}.  The mechanical Schur complement yields
  the coercive effective mass $m_{\rm eff}$, which controls $\dot p$ in
  $L^2(0,T;L^2(\Omega))$ even when $s_0=0$.  The mechanical
  reconstruction then gives the stated bounds for $\dot u$ and
  $\dot p_T$.  No estimate of $\dot p$ in $Q$ is obtained or needed.
\end{proof}

\subsection{Spatial regularity of the subsystem of displacement and total pressure}
\label{app:stokes-reg}

\begin{proposition}
  \label{prop:app:stokes-reg}
  Assume the hypotheses of Proposition~\ref{prop:app:basic-state-reg}.
  Assume additionally:
  \begin{enumerate}
  \item $\tilde f\in H^1(0,T;L^2(\Omega)^d)$.  This assumption implies
    both $\tilde f\in L^2(0,T;L^2(\Omega)^d)$ (used for the Stokes
    estimate) and that $\tilde f(0)\in L^2(\Omega)^d$ is well-defined
    as a trace (used in the compatibility condition below).  
  \item The initial data satisfy the mechanical compatibility condition:
    with $p_{T,0}:=-\lambda\operatorname{div}u_0+\alpha p_0$,
    \[
      a_u(u_0,v)+b(p_{T,0},v)
      = (\tilde f(0),v)_\Omega
      \qquad\forall v\in V.
    \]
  \end{enumerate}
  Then
  \begin{align}
    \label{eq:app:stokes-reg-est}
    u \in L^2(0,T;\boldsymbol{H}^2(\Omega)),
    \qquad
    p_T \in L^2(0,T;H^1(\Omega)).
  \end{align}
\end{proposition}

\begin{proof}
  Integrating the mechanical equation of \eqref{eq:tp:strong-pT}
  in time from $0$ to $t$ and using the mechanical compatibility
  condition (assumption~2 above) gives, for a.e.\ $t\in(0,T)$,
  \begin{align*}
    a_u(u(t),v) + b(p_T(t),v) &= (\tilde f(t),v)_\Omega
    \qquad\forall v\in V.
  \end{align*}
  From the total pressure definition \eqref{eq:total-pressure-def},
  \begin{align*}
    b(q_T,u(t)) - \lambda^{-1}(p_T(t),q_T)_\Omega
    = -\lambda^{-1}\alpha(p(t),q_T)_\Omega
    \qquad\forall q_T\in Q_T.
  \end{align*}
  The Stokes-type regularity estimate stated at the beginning of
  this appendix (assumption~(2) in the enumerated list) therefore
  applies with
  \begin{align*}
    F(t) := \tilde f(t),
    \qquad
    \zeta(t) := -\lambda^{-1}\alpha p(t).
  \end{align*}
  Since Proposition~\ref{prop:app:basic-state-reg} gives
  $p\in L^\infty(0,T;H^1(\Omega))$, we have
  $\zeta\in L^2(0,T;H^1(\Omega))$, and
  \begin{align*}
    \|u(t)\|_{\boldsymbol{H}^2} + \|p_T(t)\|_{H^1}
    \le C_{\rm St}\bigl(\|\tilde f(t)\|_{L^2}
      + \|p(t)\|_{H^1}\bigr).
  \end{align*}
  Squaring and integrating over $(0,T)$ yields \eqref{eq:app:stokes-reg-est}.
  The $H^1$-in-time parts follow from
  Proposition~\ref{prop:app:basic-state-reg}.
\end{proof}

\begin{rmk}
  The $H^1(\Omega)$-norm of $\zeta$ in assumption~(2) above is
  essential here: $\zeta(t) = -\lambda^{-1}\alpha p(t)$
  lies in $H^1(\Omega)$ because $p\in L^\infty(0,T;Q)=L^\infty(0,T;H^1_0(\Omega))$.
  If only $\zeta\in L^2(\Omega)$ were available, the Stokes estimate
  would not deliver $r\in H^1(\Omega)$.
\end{rmk}

\subsection{Spatial regularity of the diffusion block}
\label{app:elliptic-reg}

\begin{proposition}[Scalar elliptic $H^2$ regularity for $p$]
  \label{prop:app:elliptic-reg}
  Assume the hypotheses of Proposition~\ref{prop:app:basic-state-reg}.
  Then
  \begin{align}
    \label{eq:app:elliptic-reg-est}
    p \in L^2(0,T;H^2(\Omega))\cap H^1(0,T;L^2(\Omega)).
  \end{align}
\end{proposition}

\begin{proof}
  For a.e.\ $t\in(0,T)$, the state equation
  \eqref{eq:tp:strong2}
  (rewritten in total pressure form) gives
  \begin{align*}
    -\operatorname{div}(\kappa_p\nabla p(t))
    &= g_p(t) + \lambda^{-1}\alpha\dot p_T(t)
       - s_1\dot p(t)
       =: F_p(t).
  \end{align*}
  By Proposition~\ref{prop:app:basic-state-reg},
  $\dot p_T\in L^2_tQ_T$, $\dot p\in L^2_tL^2_x$,
  and $g_p\in L^2_tL^2_x$ by assumption
  \eqref{eq:app:data-energy}.  Hence $F_p\in L^2(0,T;L^2(\Omega))$.
  The scalar $H^2$-regularity estimate (assumption (1) above) gives,
  for a.e.\ $t$, $\|p(t)\|_{H^2} \le C_{\rm ell}\|F_p(t)\|_{L^2}$.
  Integrating over $(0,T)$ gives $p\in L^2(0,T;H^2(\Omega))$.
  The $H^1(0,T;L^2(\Omega))$ part was already established in
  Proposition~\ref{prop:app:basic-state-reg}.
\end{proof}

\subsection{Regularity of the optimal state}
\label{app:state-reg}

\begin{cor}
  \label{cor:app:state-reg}
  Let $\widehat m_p\in\mathcal{M}_{p,ad}$ be
  the optimal control and $(\widehat u,\widehat p_T,\widehat p)$ the
  corresponding state.  Assume \eqref{eq:app:data-energy} with
  $g_p=\widehat m_p$.  Assume also the
  hypotheses of Proposition~\ref{prop:app:stokes-reg}, namely
  $\tilde f\in H^1(0,T;L^2(\Omega)^d)$ and the mechanical
  compatibility condition
  $a_u(u_0,v)+b(p_{T,0},v)=(\tilde f(0),v)_\Omega$ for all $v\in V$,
  with $p_{T,0}:=-\lambda\operatorname{div}u_0+\alpha p_0$.
  Then
  \begin{align}
    \label{eq:app:state-reg-full}
    \widehat u &\in L^2(0,T;\boldsymbol{H}^2(\Omega))\cap H^1(0,T;V),
    \notag\\
    \widehat p_T &\in L^2(0,T;H^1(\Omega))\cap H^1(0,T;Q_T),
    \notag\\
    \widehat p
    &\in L^2(0,T;H^2(\Omega))\cap H^1(0,T;L^2(\Omega)).
  \end{align}
\end{cor}

\begin{proof}
  Combine Propositions~\ref{prop:app:basic-state-reg},
  \ref{prop:app:stokes-reg}, and \ref{prop:app:elliptic-reg}.
\end{proof}

\subsection{Regularity of the optimal adjoint}
\label{app:adjoint-reg}

\begin{proposition}
  \label{prop:app:adjoint-reg}
  Assume the state regularity of Corollary~\ref{cor:app:state-reg}
  and additionally
  \begin{align}
    \label{eq:app:data-target}
    u_C \in L^2(0,T;L^2(\Omega)^d),
    \qquad
    p_C \in L^2(0,T;L^2(\Omega)).
  \end{align}
  Then the optimal adjoint $(\widehat v^*,\widehat q_T^*,\widehat q^{\,*})$
  satisfies
  \begin{align}
    \label{eq:app:adjoint-reg-est}
    \widehat v^* &\in L^2(0,T;\boldsymbol{H}^2(\Omega))\cap H^1(0,T;V),
    \notag\\
    \widehat q_T^* &\in L^2(0,T;H^1(\Omega))\cap H^1(0,T;Q_T),
    \notag\\
    \widehat q^{\,*}
    &\in L^2(0,T;H^2(\Omega))\cap H^1(0,T;L^2(\Omega)).
  \end{align}
\end{proposition}

\begin{proof}
  Set $s:=T-t$ and define the reversed adjoint variables
  \begin{align*}
    \check v^*(s) &:= \widehat v^*(T-s),
    \quad
    \check q_T^*(s) := \widehat q_T^*(T-s),
    \quad
    \check q^*(s) := \widehat q^{\,*}(T-s).
  \end{align*}
  The terminal condition \eqref{eq:tp:adj:TC} becomes the homogeneous
  initial condition
  $(\check v^*(0),\check q_T^*(0),\check q^*(0))=(0,0,0)$,
  and the reversed adjoint system has the same coercive saddle-point
  and storage-diffusion structure as the forward state system
  (cf.\ the time-reversal argument in Section~\ref{sec:OCP:optimality}).
  The right-hand sides in $s$ are
  \begin{align*}
    \check F_v(s) &:= \alpha_{u,C}(\widehat u(T-s)-u_C(T-s))\in L^2_sL^2_x,
    \\
    \check F_q(s) &:= \alpha_{p,C}(\widehat p(T-s)-p_C(T-s))\in L^2_sL^2_x,
  \end{align*}
  by \eqref{eq:app:data-target} and the state regularity.
  Applying Proposition~\ref{prop:app:basic-state-reg} to the reversed
  system gives $\check v^*\in H^1_sV$, $\check q_T^*\in H^1_sQ_T$,
  $\check q^*\in L^\infty_s Q\cap H^1_sL^2_x$.

  For spatial regularity of $(\widehat v^*,\widehat q_T^*)$, integrating \eqref{eq:tp:adj1} from $t$ to $T$ and using the homogeneous
  terminal condition gives, for a.e.\ $t$,
  $a_u(\widehat v^*(t),v) + b(\widehat q_T^*(t),v)
  = \bigl(\int_t^T \alpha_{u,C}(\widehat u(s)-u_C(s))\,ds,\, v\bigr)_\Omega$.
  Integrating \eqref{eq:tp:adj2} from $t$ to $T$ gives
  $b(q_T,\widehat v^*(t)) - \lambda^{-1}(\widehat q_T^*(t),q_T)_\Omega
  = \lambda^{-1}\alpha(\widehat q^{\,*}(t),q_T)_\Omega$.
  The Stokes-type estimate applies with
  \begin{align*}
    F^*(t) &:= \int_t^T\alpha_{u,C}(\widehat u(s)-u_C(s))\,ds
      \in L^2(0,T;L^2(\Omega)^d),
    \\
    \zeta^*(t) &:= \lambda^{-1}\alpha\widehat q^{\,*}(t)
      \in L^2(0,T;H^1(\Omega)),
  \end{align*}
  where $\zeta^*\in L^2_tH^1_x$ because
  $\widehat q^{\,*}\in L^\infty_tQ=L^\infty_tH^1_0(\Omega)$.
  Hence $\widehat v^*\in L^2(0,T;\boldsymbol{H}^2(\Omega))$ and
  $\widehat q_T^*\in L^2(0,T;H^1(\Omega))$.

  For spatial regularity of $\widehat q^{\,*}$, the adjoint diffusion equation \eqref{eq:tp:adj3}
  can be written, for a.e.\ $t$, as a scalar elliptic equation:
  \begin{align*}
    (\kappa_p\nabla\widehat q^{\,*}(t),\nabla q)_\Omega
    &= \alpha_{p,C}(\widehat p(t)-p_C(t),q)_\Omega
      + \lambda^{-1}\alpha(\dot{\widehat q}_T^*(t),q)_\Omega
      + s_1(\dot{\widehat q}^*(t),q)_\Omega.
  \end{align*}
  The right-hand side contains
  $\widehat p-p_C\in L^2_tL^2_x$
  (by Corollary~\ref{cor:app:state-reg} and \eqref{eq:app:data-target})
  and $\dot{\widehat q}_T^*,\dot{\widehat q}^*\in L^2_tL^2_x$
  (from the energy estimate on the reversed system).
  The scalar $H^2$-regularity estimate gives $\widehat q^{\,*}\in L^2(0,T;H^2(\Omega))$.

  Since time-reversal is an isometry on all relevant Bochner spaces,
  the established regularity for the reversed variables translates
  directly to $(\widehat v^*,\widehat q_T^*,\widehat q^{\,*})$,
  yielding \eqref{eq:app:adjoint-reg-est}.
\end{proof}

\section*{Declarations}

OpenAI's ChatGPT was used to explore preliminary duality-based error estimates for the analysis in Section 5. The resulting mathematical arguments were checked and revised by the author. ChatGPT was also used to assist in writing Python scripts for the numerical experiments implemented with NGSolve. The author takes full responsibility for the mathematical analysis, numerical results, and all content of this manuscript.


\end{document}